\documentclass[10pt]{scrartcl}
\pdfoutput=1
\usepackage{ae,lmodern}
\usepackage[english]{babel}
\usepackage[utf8]{inputenc}
\usepackage[T1]{fontenc}

\usepackage[svgnames]{xcolor}
\usepackage{tikz}
\usetikzlibrary{decorations.pathreplacing, automata, arrows}
\usepackage{pgfplots, pgfplotstable}
\pgfplotsset{compat = 1.16}
\usepackage[]{graphicx}

\usepackage{mathtools}
\usepackage{amsfonts}
\usepackage{amssymb}

\usepackage{mathrsfs}
\usepackage{dsfont}
\allowdisplaybreaks

\usepackage[font=small,labelfont=bf,
labelsep=endash]{caption}
\usepackage{subcaption}
\usepackage{geometry}

\usepackage{abstract}
\usepackage{microtype}

\usepackage[bookmarks=true, breaklinks]{hyperref}
\hypersetup{
	colorlinks=true,
	linkcolor=black,
}

\usepackage{amsthm}
\newcommand{\bm}{\boldsymbol}
\hypersetup{pdfauthor={Joffrey Mathien},pdftitle={Diameter of a new model of random hyperbolic surfaces}}

\newtheorem{theorem}{Theorem}[section]
\newtheorem{proposition}[theorem]{Proposition}
\newtheorem{lemma}[theorem]{Lemma}
\newtheorem{corollary}[theorem]{Corollary}
\newtheorem{definition}[theorem]{Definition}

\newcommand{\R}{\mathds{R}}
\newcommand{\N}{\mathds{N}}
\newcommand{\Z}{\mathds{Z}}
\newcommand{\Rp}{\R^{+}}
\newcommand{\Rpn}{\R^{\ast}_{+}}

\newcommand\limit[2]{\overset{#1}{\underset{#2}{\longrightarrow}}}
\newcommand\limitRinf[1]{\limit{#1}{R \to \infty}}
\newcommand\limitginf[1]{\limit{#1}{g \to \infty}}
\newcommand\bigop[3]{\overset{#3}{\underset{#2}{#1}}}

\DeclareMathOperator{\Good}{Good}
\DeclareMathOperator{\Vol}{Vol}
\DeclareMathOperator{\supp}{Supp}
\DeclareMathOperator{\diam}{diam}
\newcommand{\p}{\mathbb{P}}
\newcommand{\E}{\mathbb{E}}

\DeclareMathOperator{\arcsinh}{arcsinh}

\title{Diameter of a new model of random hyperbolic surfaces}
\author{Joffrey Mathien\thanks{Aix Marseille Univ, CNRS, I2M, Marseille, France. email:  \href{mailto:joffrey.mathien@univ-amu.fr}{\texttt{joffrey.mathien@univ-amu.fr}}}}
\date{}

\numberwithin{equation}{section}
\numberwithin{figure}{section}

\begin{document}
	\maketitle
	
	\begin{abstract}
		The study of random surfaces, especially in the asymptotics of large genus, has been of increasing interest in recent years. Many geometrical questions have analogous formulations in the theory of random graphs with a large number of vertices, and results obtained in one area can inspire the other. In this paper, we are interested in the diameter of random surfaces, a basic measure of the connectivity of the surface. We introduce a new class of models of random surfaces built from random graphs and we compute the asymptotics of the diameter of these surfaces, which is logarithmic in the genus of the surface. The strategy of the proof relies on a detailed study of an exploration process which is the analogue of the breadth-first search exploration of a random graph. Its analysis is based on subadditive and concentration techniques.
	\end{abstract}

	\section{Introduction}
	
	An interesting quantity to study when looking at the geometry of a surface $S$ is its diameter $\diam S$. It is a basic measure of the geometry of the surface, and also provides information about objects related to the connectivity of the surface -- Cheeger's constant (\cite{Brooks90}), the first eigenvalue of the Laplacian (\cite{Cheeger, Buser82}), mixing time of Brownian motion... In this paper, we consider compact connected orientable hyperbolic surfaces without boundary, called closed hyperbolic surfaces from now on. Let us write $\mathcal M_g$ for the set of such surfaces of genus $g \geq 2$, up to isometries (called the moduli space). It is well known (see \cite{Bavard} for example) that  \[\inf_{S_g \in \mathcal{M}_g} \diam S_g \geq \ln{g} + O(1).\] In \cite{BCP_mindiam}, Budzinski, Curien, and Petri give a probabilistic proof of the sharpness of this bound, in the sense that \[\frac{1}{\ln g} \inf_{S_g \in \mathcal{M}_g} \diam S_g \limit{}{g \to \infty} 1.\] On the other hand, it is possible to build hyperbolic surfaces of any genus $g \geq 2$ with arbitrarily large diameter using the collar lemma (see Appendix \ref{collar} for the lemma and an idea of the construction of such a surface).
	
	A natural question is therefore the "typical" behaviour of the diameter of a surface: when choosing a surface "at random", what does its diameter look like? In particular, is, like for graphs (see \cite{Bollobas}), the typical behaviour close to the optimal one? To answer this question, different models of random surfaces have already been considered.
	
	In the Brooks-Makover model (introduced in \cite{Brooks}), one considers surfaces built by compactifying a random gluing of ideal triangles. The resulting surfaces $\left(S_g\right)_g$ are known to be Bely\u{\i} surfaces whose genus is approximately $g$ (see \cite{Gamburd}). It is well known that Bely\u{\i} surfaces of genus $g$ form a dense subset of $\mathcal M_g$ (\cite{Belyi}). In \cite{BCP_Bel}, it is proven that for a sequence of random surfaces $S_g$ obtained with the Brooks-Makeover model, we have \[\frac{\diam S_g}{\ln g} \limit{}{g \to \infty} 2\] in probability.
	
	Another well-studied model, the so-called Weil-Petersson model, is based on the Weil-Petersson metric on the moduli space $\mathcal M_g$ (\cite{Weil_58}, \cite{WMF}). It can be related to the pair-of-pants decomposition of surfaces (see \cite[Chapters 1 and 6]{Buser}, \cite[Chapter 3]{MonkPHD}, or the description of our model below), and is therefore particularly geometric. In addition, contrary to the Brooks-Makover model, any surface $S_g \in \mathcal M_g$ can be obtained with this model. Since the work of Mirzakhani (\cite{MIR15, Mirzakhani_07b}) which provides probabilistic tools to work with it, the Weil-Petersson model has been of growing interest. In particular, in \cite[Theorem 4.1]{Mirzakhani7}, Mirzakhani proves that asymptotically almost surely (\textit{i.e.} with probability converging toward 1 as $g \to \infty$), \[\inf_{S_g \in \mathcal{M}_g} \diam S_g \leq 40 \ln g.\]
	This bound has since been improved using spectral estimates (see \cite{Magee, AM_spectralgap}). However, to the author's knowledge, there is no proof of a convergence, \[\frac{\diam S_g}{\ln g} \limit{}{g \to \infty} C\] for some deterministic $C>0$.
	
	All these results converge towards the fact that "typically" a surface has a diameter of order $\ln g$. In this paper, we study an extension of the model of surfaces introduced in \cite{BCP_mindiam}, inspired by the construction in terms of pants decomposition of surfaces, as for the Weil-Petersson model. This model, compared with the Weil-Petersson one, has the advantage to be easier to handle: it carries a more manipulable measure, making calculations easier. Nevertheless, it keeps the possibility to theoretically generate any surface in $\mathcal M_g$.
	
	\subsection{Model and Results}
	
	\begin{figure}[ht]
		\centering
		\includegraphics[width=\textwidth]{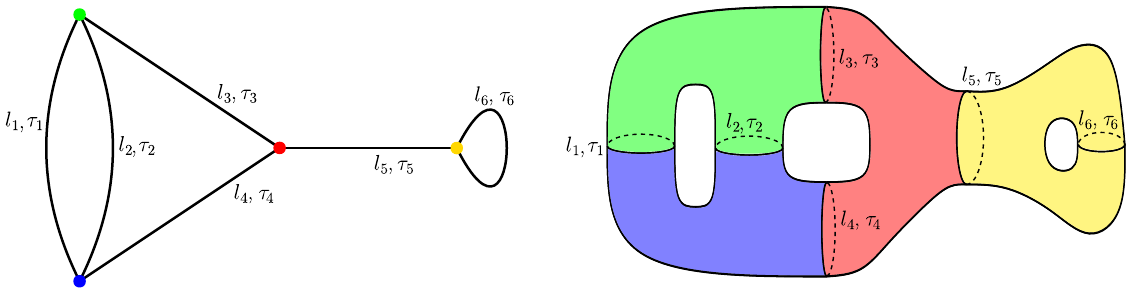}
		\caption{The correspondence between graphs and surfaces built with pair of pants.}
		\label{fig:Surfacegraph}
	\end{figure}
	
	Given a 3-regular (multi)-graph $G$ on $2g-2$ vertices, one can build a closed topological surface of genus $g$ by gluing $2g-2$ pairs of pants according to this graph: each vertex $v$ corresponds to a pair of pants $P(v)$, and an edge between two vertices $u$ and $v$ to a gluing along a boundary component between $P(u)$ and $P(v)$. The gluing is done in such a way that the surface obtained is orientable. By prescribing additional parameters, it is possible to uniquely describe a hyperbolic metric on such a surface: \begin{itemize}
		\item giving the lengths $l \in \Rpn$ of the boundary components of the pairs of pants, we describe a unique hyperbolic metric on each pair of pants, with geodesic boundary components.
		\item giving a twist parameter $\tau \in \R/2\pi\Z$ for each gluing -- which corresponds to the choice of a parameterization of the boundaries we glue -- we describe how to "glue" the metrics on pairs of pants together to get a global hyperbolic metric.
	\end{itemize} 
	This gives $6g-6$ parameters, called Fenchel-Nielsen coordinates: one length and one twist for each gluing or equivalently for each edge in the graph. Therefore, it can be summarized in an edge-weighted graph, where the weight $\left(l_e, \tau_e\right)$ of an edge $e$ corresponds to the couple of parameters (length, twist) associated with this edge (see, for example, Figure \ref{fig:Surfacegraph}).
	The question of knowing precisely what it means to have twist equal to 0, \textit{i.e.} the choice of the parameterization of the boundary components of the pairs of pants, is postponed to Section \ref{behavalpha}, which is the only part in which we really need to make it explicit. More details can be found in \cite[Sections 1.7, 6.2]{Buser}.
	
	\begin{figure}[ht]
		\centering
		\includegraphics[width=0.8\textwidth]{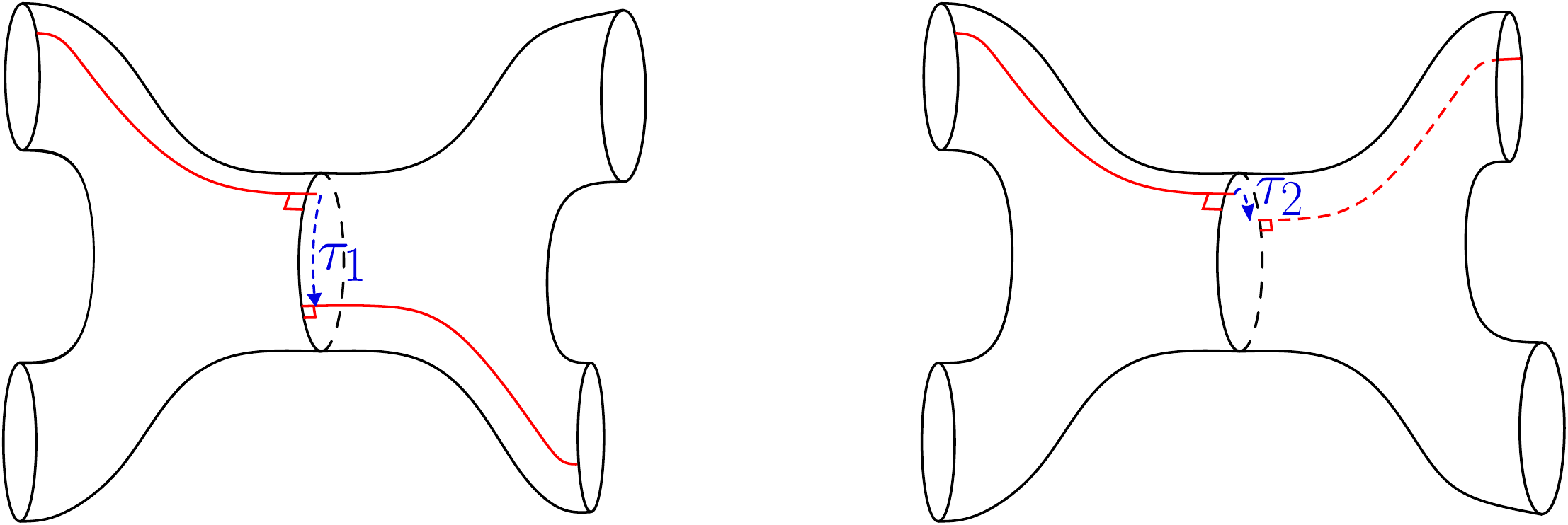}
		\caption{Two pairs of pants glued together with two different twists.}
		\label{fig:twist}
	\end{figure}
	
\begin{figure}[ht]
		\centering
		\includegraphics[width=\textwidth]{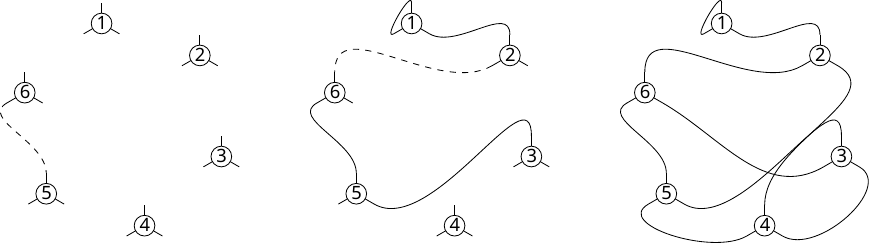}
		\caption{An example of a 3-regular graph with 6 vertices obtained with the configuration model.}
		\label{fig:config}
	\end{figure}

	For the Weil-Petersson model, the law of twists and lengths is not totally explicit -- it corresponds to the Lebesgue measure on a set that cannot be explicitly described in general -- which makes computations complicated. In \cite{BCP_mindiam}, they take the underlying graph to be random according to the configuration model (see Definition \ref{defi:config}), and fix twists to zero, and lengths to some $a >0$. We recall the definition of the configuration model for a 3-regular multigraph with $2n$ vertices:
	
	\begin{definition}\label{defi:config}
		To build a 3-regular random (multi-)graph with $2n$ vertices according to the \textbf{configuration model}, start with the $2n$ vertices, each of which has $3$ half-edges attached. Let $\mathcal{H}$ be the set of half-edges obtained.
		Choose a pairing of the half-edges uniformly at random, that is a partition of $\mathcal H$ into pairs. These pairs form the edges of the multigraph.
		We can proceed algorithmically as follows:
		\begin{itemize}
			\item  Take arbitrarily an unpaired half-edge $x \in \mathcal{H}$.
			\item Draw some half-edge $\eta(x)$ uniformly at random from all the other unpaired half-edges.
			\item The pair $\left(x, \eta(x)\right)$ then forms an edge in the graph.
			\item Continue until no half-edge remains unpaired.
		\end{itemize}
	\end{definition}
	An example is given in Figure \ref{fig:config}.
	Note that this construction attaches no importance to the order in which the edges are paired. In fact, all configurations, \textit{i.e.} all possible pairings of half-edges, are equiprobable with this construction.
	For more details on the configuration model, see \cite[Chapter 7]{hofstad_2016}, or \cite{Wormald_1999}. 
	
	In this article, we extend the model of \cite{BCP_mindiam} as follows:
	\begin{definition}\label{defi:mod}
		Given $g \geq 2$ and $\nu = \nu_g$ a probability distribution on $\Rpn\times [0, 2\pi[$, let $\mathfrak S_{\nu, g}$ be the surface built as explained before, with the underlying graph being a 3-regular graph obtained with the configuration model, and with all the gluings being independent of the graph, iid of law $\nu$.
	\end{definition}
	
	Compared with \cite{BCP_mindiam}, randomizing twists and lengths allows us to work on a richer class of surfaces, but also makes some geometric interpretations crucial in \cite{BCP_mindiam} more difficult. 
	
	The paper is mainly devoted to proving the following theorem.
	\begin{theorem}\label{th:Principal3}
		Let $\nu_l$ be the length distribution. Suppose that \[\supp\left(\nu_l\right) \subset [a_g^-, a_g^+], \quad \text{with } \frac{1}{a_g^-} = o(\ln \ln g),\  a_g^+ = o (\ln \ln g).\] Then, there exists $0< \alpha_{\nu} \leq 1$ which depends only on the law $\nu$ of the (length, twist)-weights of the pairs of pants, such that for all $\eta>0$,  the sequence \[\frac{\diam \mathfrak S_{\nu, g} - \frac{1}{\alpha_{\nu}}\ln g}{\ln^{3/4+\eta}g} \limitginf{} 0\] in probability.
		In addition, $\alpha_\nu \geq \frac{\ln 2}{\Delta_+}$, where $\Delta_+$ is any upper bound for the diameter of all pairs of pants with boundary length in $[a_g^-, a_g^+]$.
		
		Moreover, when $\supp{\nu}$ is uniformly bounded in $g$, we can replace $\ln^{3/4+\eta} g$ by $u_g \ln^{3/4} g $ where $(u_g)_g$ is any sequence such that $u_g\limitginf{} \infty$. 
	\end{theorem}
	In particular, we obtain that for any $\eta>0$,  asymptotically almost surely, $\diam \mathfrak S_{\nu, g} = \frac{1}{\alpha_\nu} \ln g + o\left(\ln^{3/4+\eta} g\right)$.
	In addition, $\alpha_\nu$ can be geometrically interpreted as the growth rate of a kind of tree-like surface, which will be introduced in Section \ref{TLS}. In Subsection \ref{behavalpha}, we partially deal with the behaviour of $\alpha$ by showing the following theorem. \begin{theorem}\label{th:balpha}
		Let us denote $\alpha_{l\otimes\nu_t}$ the value of $\alpha_\nu$ when the length distribution is the Dirac mass at $l$ and $\nu_t$ is the twist distribution.
		Then, if $\nu_t$ is uniform, \[\alpha_{l\otimes\nu_t} \limit{}{l \to \infty} 1.\]
	\end{theorem}
	
	\subsection{Idea of the proofs} The proof of Theorem \ref{th:Principal3} is mainly inspired by what Bollobás and Fernandez de la Vega did in \cite{Bollobas} to study the diameter of random graphs, and which has already been reused in \cite{BCP_mindiam}. A kind of birthday paradox argument tells us that until any two increasing balls contain $\sqrt{g}$ pairs of pants, they are not likely to intersect. On the contrary, as soon as they contain a bit more than $\sqrt{g}$ pairs of pants, any two such balls have a high probability to intersect. As a consequence, their common radius should be close to half the diameter of the surface. Then the issue is to estimate the volume growth of metric balls in the surface or, equivalently, of the number of pairs of pants in such a ball. 
	
	As a graph in the configuration model "looks like" a 3-regular tree at local scale -- when looking at sufficiently small balls -- we show in Section \ref{Surface} that $\mathfrak S_{\nu, g}$ mainly "locally behaves like" a kind of "tree-like surface" $\widehat{\mathfrak S}_\nu$, which is obtained by gluing along an infinite 3-regular tree pairs of pants with random boundary lengths and twists of joint law $\nu$. The balls of $\mathfrak S_{\nu, g}$ and $\widehat{\mathfrak S}_\nu$ are very similar, and then their growth rate is roughly the same, leading to similar diameters. As a consequence, in Section \ref{TLS}, we introduce properly $\widehat{\mathfrak S}_\nu$ and study its exponential growth rate. In particular, we show that the number of pairs of pants intersecting a ball of hyperbolic radius $R$ concentrates around $e^{\alpha_\nu R}$ for some deterministic $\alpha_\nu$. It mainly rests on concentration techniques and on the tree structure of $\widehat{\mathfrak S}_\nu$, which allows us to use subadditive results.
	
	As we will see in Section \ref{Surface}, contrary to \cite{Bollobas} and \cite{BCP_mindiam}, it is not straightforward to get a lower bound on the diameter by a direct comparison with $\widehat{\mathfrak S}_\nu$, because it is not an isometric cover of $\mathfrak S_{\nu, g}$. Similarly to the upper bound, a detailed study of a breadth-first exploration of the graph is needed.
	
	As for the proof of Theorem \ref{th:balpha}, it amounts to showing that from a metric point of view, large enough balls in the tree-like surface behave as in the hyperbolic plane. In particular, the tree structure and a large injectivity radius allow us to describe easily a set of pairs of pants around which the balls have a sufficient growth rate. The uniformity of the twists is used to ensure that there are enough of these pairs of pants to have a significant effect on the global growth of the surface.
	
	\subsection{Remarks}
	Note the fluctuations of order $\ln^{3/4} g$ obtained for the diameter. These come from the concentration bound of Theorem $\ref{th:concentration}$ and are not optimal. In particular, if we work harder in the proof of Lemma \ref{lem:Majca}, we could end up with fluctuations of order $\ln^{1/2} g$ (see the remark at the end of the lemma). Even in this case, they are not expected to be optimal: in fact, we believe that the good order for fluctuations should be roughly of the order $\ln \ln g$, as for random graphs (see \cite{Bollobas}).
	
	Note also the limits on the support of $\nu_l$, and more especially the upper bound: to generate any surface, it is not necessary to have a support equal to $\Rp$. Indeed, there exists a constant $B_g$, called Bers' constant, such that any surface of genus $g$ can be decomposed into pairs of pants with all boundary components of length less than $B_g$. However, $B_g$ is known to be at least of order $\sqrt{g}$ (\cite[Theorem 5.1.3]{Buser}) – and at most of order $g$ (see e.g.  \cite{Parlier} for the best known bounds) -- and thus, our theorem fails to reach this bound. If we were able to improve the concentration inequality of Theorem \ref{th:concentration}, and especially the minimal radius $R_0$ of the sphere for which it holds, we could obtain better bounds. In particular, remarks after Lemmas \ref{lem:Majca} and \ref{lem:cMDiarmid} go in this direction. Nevertheless, we do not expect to reach $B_g$ by this method. Indeed, the fact that we are using subadditive techniques should limit the proof to pairs of pants whose boundary lengths are small compared to the diameter of the surface, \textit{i.e.} $\ln g$. Heuristically, subadditivity tells us that a large-scale phenomenon (the diameter) is obtained thanks to local contributions (whose scale is the boundary length of the pairs of pants).
	
	For the sake of clarity, the author decided to state a rather general statement. By following the proof with more information on the law $\nu$, slightly more precise estimates on the fluctuations could be obtained.
	
	An interesting question that is barely answered in this paper is the behaviour of $\alpha_{\nu}$. In particular, we prove that $\alpha_{l \otimes \nu_t} \to 1$ as $l \to \infty$ when twists are uniform and $l$ is the deterministic length of the boundary components of the pairs of pants, and conjecture that this should be true for many twist laws. On the contrary, we give no information about what happens for random lengths in general. In particular, the question of the continuity of $\alpha_\nu$ with respect to $\nu$ is not answered in the paper.

	\subsection{About the notations}
	Throughout the article, and for ease of reading, we try to use bold fonts for the notations the \textbf{first} time we introduce them, to make it easier to find the definition. Then we do not bold them again.
	
	Consider a graph $G$. In what follows, we use the following conventions and notations: \begin{itemize}
		\item We use \textbf{the exponent $\bm{G}$} to refer to any object "linked" with $G$: for example, $V^G$ will denote the set of vertices of $G$. When no confusion is possible, we will omit it.
		\item For any surface $S^G$ built from a graph $G$, we make no distinction between a vertex $v$ of $G$ and the corresponding pair of pants $P(v)$ used to build $S^G$, using the point of view "surface $S^G \approx \text{ edge-weighted graph } (G = (V, E), w = (l_e, t_e)_{e \in E})$". In particular, we often write $P$ indistinctly for $v$ and $P(v)$.
		\item As a consequence, we use the notation $\bm{d_h^G}$ for the distance induced by the hyperbolic metric on the surface $(G, w)$.
		\item More generally, we allow ourselves to assimilate any notion on $V^G$ and on the set of pairs of pants used to build $S^G$: for example, if $G$ is a rooted tree, we call \textbf{height} of the pair of pants $P$ the height in the tree of the associated vertex. Similarly, we allow ourselves to write $x \in v$ to mean $x \in P(v)$.
	\end{itemize}
	
	We recall also that \begin{itemize}
		\item for $S$ a finite set, $\bm{\#S}$ denotes its cardinality.
		\item for any measurable set $A$, its indicator function is denoted $\bm{{\mathds{1}_{A}}}$.
	\end{itemize}

	\subsection{Acknowledgment}
	The author would like to thank Adrien Boulanger and his PhD advisor Charles Bordenave for their involvement and support, and for our fruitful discussions. He also thanks the referee of the paper for many useful comments.
	
	\section{Tree-like surface} 
	\label{TLS}
	
	On a small scale, graphs obtained with the configuration model look like trees. Therefore, we prove in Section \ref{Surface} that the surfaces $\mathfrak S_{\nu, g}$ we work on locally behave, in some sense, like a surface built from a tree. In this section, the law $\nu$ is fixed and we analyze this "tree-like surface" to prepare the study of $\mathfrak S_{\nu, g}$ in Section \ref{Surface}.
	
	\subsection{Notations and main results}
\begin{figure}[ht]
		\centering
		\includegraphics[width=0.45\textwidth]{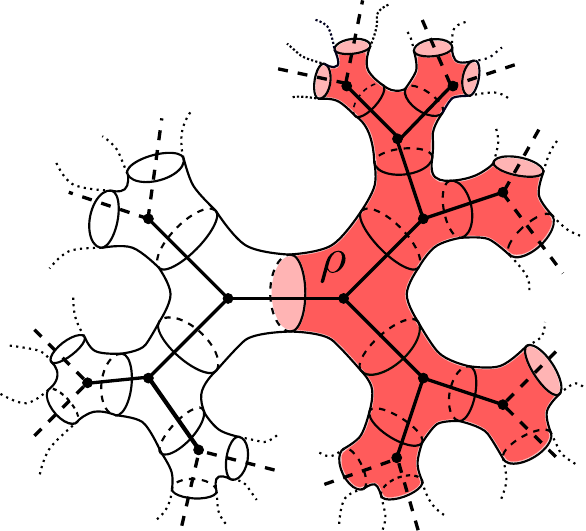}
		\caption{The surface $\widehat{\mathfrak S}_\nu$, build from $\mathfrak T_3$, and in red, the surface $\mathfrak S_\nu$, obtained by only considering $\mathfrak B$.}
		\label{fig:def1}
	\end{figure}

	Let $\bm{{(\mathfrak T_3, \rho)}}$ denote the infinite 3-regular tree rooted at $\rho$, and $\bm{{(\mathfrak B, \rho)}}$ the infinite binary tree rooted at $\rho$. Let $\nu$ be a probability distribution on $\Rpn\times [0, 2\pi[$.  In this section, we study the two following surfaces (see also Figure \ref{fig:def1}).
	
	\begin{definition}
		\begin{itemize}
			\item Let $\bm{{\widehat{\mathfrak S}_\nu}}$ be the surface obtained by gluing pairs of pants according to $\mathfrak T_3$, with lengths and twists  following the law $\nu$, independently from one edge to another.
			\item Let $\bm{{\mathfrak S_\nu}}$ be the surface obtained by gluing pairs of pants according to $\mathfrak B$, with lengths and twists following the law $\nu$, independently from one edge to another.
		\end{itemize}
	\end{definition}
	
	For the whole section, we focus on the case when the length distribution is compactly supported in $\left [2l^-, 2l^+\right]$, and try to make the dependence on $l^-$ and $l^+$ explicit.
	
	\begin{definition}
		Let $P$ be a pair of pants. We denote :\begin{itemize}
			\item $\bm{{\delta_P}}$ the minimal distance between two boundary components of $P$.
			\item $\bm{{\Delta_P}}$ the diameter of $P$.
		\end{itemize}
	\end{definition}
	These two geometrical quantities will be of great interest for our study, to approximate the growth of the surface by the one of the underlying tree. In particular, with our assumptions on the length distribution, we can obtain uniform bounds on $\delta$ and $\Delta$ which only depend on $l_-$ and $l_+$.
	
	\begin{proposition}\label{prop:Ddelta}
		{Suppose that a pair of pants P has its boundary lengths bounded between $2l^-$ and $2l^+$. Then: \begin{enumerate}
				\item the minimal distance $\delta_P$ between two boundary components of $P$ is bounded in \[\left[ \delta_- = \cosh^{-1}\left(\frac{\cosh l^- + \cosh^2 l^+}{\sinh ^2 l^+} \right),\ \delta_+ = \cosh^{-1}\left(\frac{\cosh l^+ + \cosh^2 l^-}{\sinh ^2 l^-} \right)\right].\]
				\item In addition, \[\delta_-\geq \exp(-2l^+) \quad \text{and}\quad \delta_+ \leq 2l^+ +\ln(4)-2\ln l^-.\]
				\item The diameter $\Delta_P$ of $P$ is bounded from above by some $\Delta_+$ such that: \[\max{\left({l^+}, \delta_+\right)} \leq \Delta_+ \leq 8\max(l^+, \delta_+).\]
				\item In addition, \[\Delta_+ \geq 1.31 \quad \text{and}\quad \ln\left(\frac{\Delta_+}{\delta_-} +1\right) \leq 4\Delta_+.\]
		\end{enumerate}}
	\end{proposition}
	
	The proof of the proposition, which is only based on hyperbolic trigonometry on the hyperbolic plane, is postponed to Appendix \ref{pantalon}.
	
\begin{figure}[ht]
		\centering
		\includegraphics[width=0.9\textwidth]{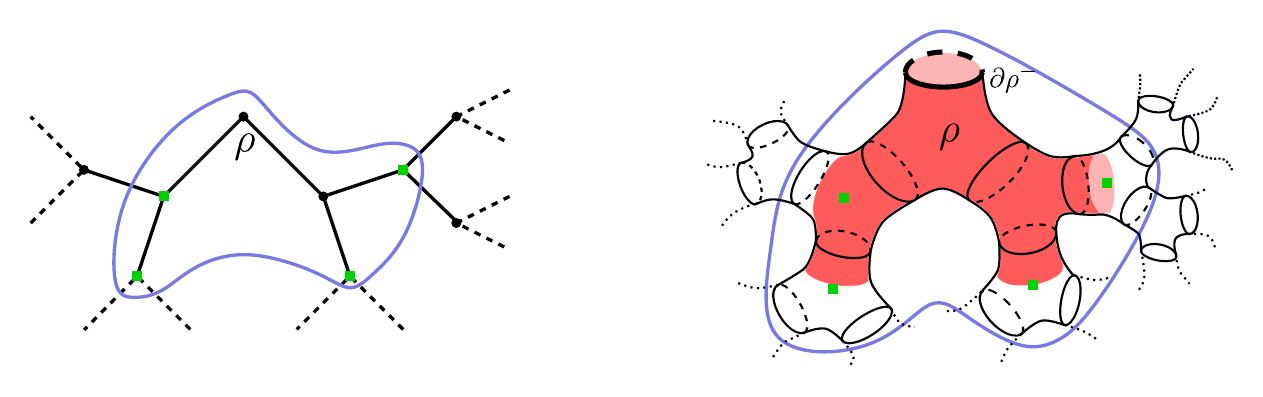}
		\caption{The region coloured on the right corresponds to the points at a distance less than $R$ from $\partial \rho^-$. The ball $\mathcal{B}_R$ is surrounded in blue. The squared vertices, in green, correspond to the pairs of pants in $\mathcal S_R$. In this example, $N_R = 4$.}
		\label{fig:sphere}
	\end{figure}

	We want to understand the growth rate of $\widehat{\mathfrak S}_\nu$, that is, roughly, the number of pairs of pants in a hyperbolic ball or sphere. In order to simplify the problem (at least for the proof we give here), we are going to estimate the growth in only two of the three branches of the tree-like surface. Up to a factor of $\frac{3}{2}$, this gives the correct growth and does not change the exponential rate (see later Figure \ref{fig:treelike_sum}, and Corollary \ref{coro:surfacecomplete}). As a consequence, we focus on $\mathfrak S_\nu$. We are interested in the following quantity (see also Figure \ref{fig:sphere}):
	\begin{definition}
		Let $\bm{{\partial \rho^{-}}} = \partial \mathfrak S_\nu \subset \partial \rho$ for the unglued boundary component of $\rho$, which is also the boundary of $\mathfrak S_\nu$, and let $R \in \Rp$. We define \[\bm{{\mathcal B_R}} = \left\lbrace P \in V^{\mathfrak B} \text{ such that $d_h(\partial \rho^-, P) \leq R$}\right \rbrace, \] \[\bm{{\mathcal S_R}} = \left\lbrace P \in \mathcal B_R \text{ such that there exists $\widetilde{P}$ a child of $P$, $\widetilde{P} \notin \mathcal B_R$}\right\rbrace,\] and\[\bm{{N_r}} = \#\mathcal S_R.\]
	\end{definition}
	
	As we will see later on, we should think of $\mathcal B_R$ (respectively $\mathcal S_R$) as a kind of ball (resp. sphere) of radius $R$ centered at $\rho$. The aim of the section is to understand the growth of $N_R$ as a function of $R$, and more precisely to prove the following result.
	\begin{theorem}[Subgaussian concentration]\label{th:concentration} Suppose that the law of the boundary lengths is compactly supported:	\[\supp\left(\nu_l\right) \subset [2l^-, 2l^+].\]
		Then there exists some $\alpha = \alpha_{\nu}>0$,  which depends only on $\nu$, and a numerical constant $C >0$ independent of all other variables such that
		\[\forall 0< \varepsilon <1, \exists R_0 > 0, \forall R > R_0, \quad \p\left(\left\lvert \frac{\ln N_R}{R} - \alpha\right\rvert \geq \varepsilon\right) \leq 2 \exp\left(-e^{-C\Delta_+}\varepsilon^2R \right).\]
		In addition, there exists a numerical constant $c >0$ independent of all other variables such that \[R_0 \leq \frac{e^{c\Delta_+}}{\varepsilon^4}.\]
	\end{theorem}
	
	Note that this implies that \[\label{IneqThmC}\tag{\ref{th:concentration}.1} \forall 0 < \varepsilon <1, \forall R \geq R_0, \quad 2 \exp\left(-e^{-C\Delta_+} \varepsilon^2R \right) \geq \begin{cases}\p\left(N_R \leq e^{(\alpha-\varepsilon) R}\right)\\ \p\left(N_R \geq e^{(\alpha+\varepsilon) R}\right) \end{cases} .\] We will mainly use this last form of the concentration in Section \ref{Surface}.
	
	To prove Theorem \ref{th:concentration}, we need to introduce additional definitions and notations which will be used later on in the section. In particular, we extend the concepts of sphere and of $N_R$ to any pair of pants in $V^\mathfrak B$ (not only for $\rho$).
	First, we recall that by abuse, we generally assimilate a pair of pants and the corresponding vertex in the graph used to construct the surface considered.
	
\begin{figure}[ht]
		\centering
		\includegraphics[width=0.8\textwidth]{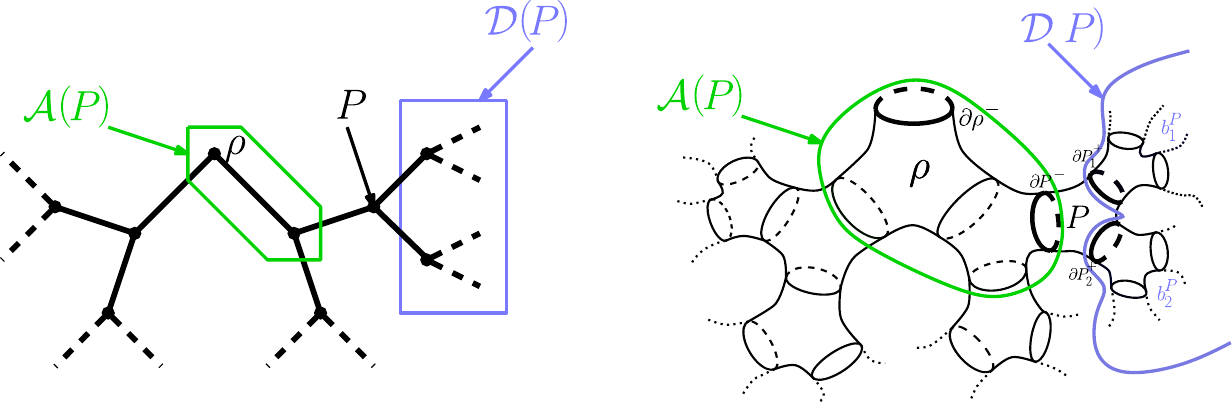}
		\caption{$\mathcal A(P)$ is surrounded by green, $\mathcal D(P)$ by blue.
			The boundary components $\partial \rho^-$, $\partial P^-$ and $\partial P_i^+, i=1, 2$ are shown in bold on the right.}
		\label{fig:surfacearborescente}
	\end{figure}

	\begin{definition} (See Figure \ref{fig:surfacearborescente})
		\begin{itemize}
			\item We define a partial \textbf{order relation} on $V^\mathfrak B$, for which $\rho$ is a maximum: \[\bm{u < v} \text{ if $u$ is a descendant of $v$}.\]
			\item We will use these ancestry/descendants notions for pairs of pants and so we note \[\bm{{\mathcal A(P)}} = \left \lbrace \widetilde{P} \in V^\mathfrak B, \widetilde{P} > P \right \rbrace, \qquad \bm{{\mathcal D(P)}} = \left \lbrace \widetilde{P}\in V^\mathfrak B, \widetilde{P} < P \right \rbrace 
			\]and \[\bm{{\overline{\mathcal A}(P)}} = \left \lbrace \widetilde{P} \in V^\mathfrak B, \widetilde{P} \geq P \right \rbrace = \mathcal A(P) \cup P,
			\qquad \bm{{\overline{\mathcal D}(P)}} = \left \lbrace \widetilde{P}\in V^\mathfrak B, \widetilde{P} \leq P \right \rbrace = \mathcal D(P) \cup P.\]
			\item More generally, if $S$ is a set of pairs of pants, we write \[\bm{{\overline{\mathcal A}(S)}} = \bigop{\bigcup}{P \in S}{} \overline{\mathcal A}(P).\]
		\end{itemize}  
		
	\end{definition}
	
	\begin{definition}
		(See Figure \ref{fig:surfacearborescente}) In the following, \begin{itemize}
			\item for a pair of pants $P \neq \rho$, we denote $\bm{{{\partial P^-}}}$ the set of points which are both in $P$ and $\mathcal A(P)$.
			\item for any $P$, we denote $(\bm{{\partial P^{+}_{i}}},\ i = 1, 2)$ the two connected components of the intersection of $P$ and $\mathcal D(P)$. The assignment of the labels 1 and 2 is arbitrary but do not depend on the randomness of the surface.
		\end{itemize}
	\end{definition}
	
	\begin{definition}\label{defi:ball}
		Let $P \in V^\mathfrak B$. \begin{itemize} \item We call \textbf{ball of radius $R$ and center $P$}, denoted $\bm{{\mathcal B_R(P)}}$, the set 
			\[\mathcal B_R(P) = \left \lbrace \widetilde{P} \in \overline{\mathcal{D}}(P), \  d_h(\partial P^{-}, \partial \widetilde{P}^{-}) \leq R \right \rbrace.\]
			This corresponds to the pairs of pants in the descendancy of $P$ with at least one point at distance less than $R$ from $\partial P^-$. 
			\item We call \textbf{sphere of radius $R$ and center $P$}, the subset  \[\bm{{\mathcal S_R(P)}} = \left \lbrace \widetilde{P} \in \mathcal B_R(P), \  \max\left(d_h(\partial P^{-}, \partial \widetilde{P}^{+}_1),\ d_h(\partial P^{-}, \partial \widetilde{P}^{+}_2)\right)> R \right\rbrace.\]
			These are the pairs of pants of $\mathcal{B}_R(P)$ such that at least one of their sons is not in $\mathcal{B}_R(P)$, and therefore, which are "roots" of an unseen branch of the surface. It is in some sense the boundary of the ball $\mathcal B_R(P)$ (see Figure \ref{fig:sphere}).
			\item The cardinality of $\mathcal S_R(P)$ is denoted  $\bm{{N_R(P)}}$. %When $R<0$, by convention $N_R(P) = 1$.
			\item Following the usual terminology for partially ordered sets, we say that a family $(P)_{P\in I}$ of pairs of pants is an \textbf{anti-chain} if there is no ancestry/descendancy relationship between the elements of $I$ (\textit{i.e.} for any $P, \widetilde{P} \in I$, $\widetilde P \notin \overline{\mathcal D}(P)$ and $P \notin \overline{\mathcal D}(\widetilde{P})$).
			In such a situation, the family $(N_R(P))_{P \in I}$ is a family of iid random variables by construction.
		\end{itemize}
	\end{definition}
	
	In particular, $\mathcal B_R = \mathcal B_R(\rho)$, $\mathcal S_R = \mathcal S_R(\rho)$ and $N_R = N_R(\rho)$.
	
	We finally introduce two last sets of pairs of pants, which will be helpful in what follows. They both "surround" the set $\mathcal S_r$ with the convenient property that when looking at all their ancestor (see Definition \ref{defi:treeabove}), we get a complete subtree in the sense that a node is either a leaf or a saturated node. It will be convenient to state a Markov property later.
	
	\begin{definition}\label{defi:U}
		For $R \geq 0$,
		\begin{itemize}
			\item let $\bm{{U_R}}$ be the set of maximal pairs of pants of $\mathcal S_R$, \textit{i.e.} the subset of pairs of pants in $\mathcal S_R$ with no ancestor in $\mathcal S_R$. In particular, $U_R$ is an anti-chain;
			\item let \[\bm{{U'_R}} = \left\lbrace P \in V^\mathfrak B, \exists \widetilde P \in \mathcal S_R,\ \text{$P$ is a child of $\widetilde P$}\right\rbrace,\] the children of elements in $\mathcal S_R$.
		\end{itemize}
		
	\end{definition}
	
	The proof of Theorem \ref{th:concentration} will be divided into two parts. First, in Section \ref{convexp}, we show that the exponential growth of $N_R$ converges in expectation to some $\alpha$  thanks to subadditive techniques based on the tree structure. Then, in Section \ref{convas}, we use a concentration inequality to extend the result to $N_R$ itself. Section \ref{behavalpha} is devoted to the proof of Theorem \ref{th:balpha}.

	\subsection{Convergence in expectation}\label{convexp}
	
	The aim of this subsection is to prove the following property and corollary, which will be key points for Theorem \ref{th:concentration}.
	
	\begin{proposition}[Submultiplicativity of the process]\label{prop:ssmult}
		Fix $R \geq 0$ and recall $U_R$ from Definition \ref{defi:U}.
		Then \[\forall r \geq 0,\quad  N_{R+r} \leq 2e^{4\Delta_+} \bigop{\sum}{P \in U_{R}}{} N_{r}(P).\label{ssmult}\tag{\ref{prop:ssmult}.1}\]
	\end{proposition}
	
	\begin{corollary}\label{coro:convexp}
		There exists $\alpha = \alpha_\nu > 0$ such that \[\Sigma_R = \frac{\ln \E N_R(\rho)}{R} \limitRinf{}{} \alpha,\quad  \text{with }\left\lvert \alpha - \Sigma_R \right\rvert \leq \frac{15\Delta_+}{R}. \label{ratealpha}\tag{\ref{coro:convexp}.1} \]
	\end{corollary}
	
	The subsection is divided as follows: first, we give some general facts about the objects we defined in the previous subsection, which are mainly based on some simple geometrical or combinatorial arguments. Then we prove some deterministic bounds from these properties and the tree structure of $\mathfrak S_\nu$. Finally, we use these properties to prove probabilistic results on $N_R$, which leads to Property \ref{prop:ssmult} and Corollary \ref{coro:convexp}.
	
	\subsubsection{Basic facts about \texorpdfstring{$\mathfrak S_\nu$}{S\_ν}}
	
	In this section, we state some basic properties of the objects defined in the previous section that will be useful several times in the rest of the section.
	
	\begin{proposition}\label{prop:elem}
		Let $P, \widetilde P \in V^\mathfrak B$, $\widetilde {P} \in \mathcal D(P)$. Then:
		\begin{enumerate}
			\item \label{elem:dist} $d_h(\partial P^-, \widetilde{P}) = d_h(\partial P^-, \partial \widetilde{P}^-)$.
			\item \label{elem:delta} $d_h(\partial \rho^-, P) \geq R \Rightarrow d_h(\partial \rho^-, \widetilde{P}) \geq R + \delta_-$.
			\item \label{elem:Delta} If moreover, $\widetilde{P}$ is a child of $P$, then $d_h(\partial \rho^-, P) \leq R \Rightarrow d_h(\partial \rho^-, \widetilde{P}) \leq R + \Delta_+$.
			\item \label{elem:S} If $P \in \mathcal S_R$, then $d_h(\partial \rho^-, P) > R-\Delta_+$.
		\end{enumerate}
	\end{proposition}
	
	\begin{proof}
		Because $\mathfrak S_\nu$ is built from a tree,  any path from $P$ to $\widetilde{P}$ must intersect $\partial \widetilde{P}^-$, leading to the first point.
		Take now a path from $\partial \rho^-$ to $\widetilde{P}$. It induces a path on the tree $\mathfrak B$, crossing $P$. So it can be decomposed in a path from  $\partial \rho^-$ to $\partial P^-$, of length at least $R$, a path from $\partial P^-$ to one of the $\partial P_i^+$, of length at least $\delta_-$, and a last part from this $\partial P_i^+$ to $\widetilde{P}$. For the third point, take a minimal path from $\partial \rho^-$ to $\partial P^-$. The endpoint of this path is at a distance at most $\Delta_+$ of $\partial P_i^+$ for $i =1, 2$, \textit{i.e.} of the sons of $P$. Hence, an extension of the first path in a path from $\partial \rho^-$ to $\partial \widetilde{P}^-$ of length at most $R+\Delta_+$. As a consequence, if $d_h(\partial \rho^-, P)\leq R-\Delta_+$, then both sons of $P$ are in $\mathcal B_R$, so $P \notin S_R$. 
	\end{proof}
	
	Let us also make more explicit the relation between the metric structures on $\mathfrak B$ and on $\mathfrak S_\nu$.
	
	\begin{proposition}\label{prop:quasiiso}
		Let $\mathrm{Proj}: \mathfrak S_\nu \to \mathfrak B$ be the map sending a point of $ \mathfrak S_\nu$ to the nearest pair of pants from $\rho$ in $\mathfrak B$ it belongs to. Then \[\forall x, y \in  \mathfrak S_\nu, \quad \Delta_+^{-1}d(x, y)-1 \leq d(\mathrm{Proj}(x), \mathrm{Proj}(y)) \leq \delta_-^{-1}d(x, y)+1.\] In particular, $\mathrm{Proj}$ is a $\left(\Delta_+ + \delta_-^{-1},1\right)$-quasi-isometry.
	\end{proposition}
	
	\begin{proof}
		Take $x$ and $y$ in $ \mathfrak S_\nu$, and $u =\mathrm{Proj}(x), v = \mathrm{Proj}(y)$. Consider $v_0 = u, v_1, v_2 \ldots, v_{k-1}, v_{k} = v$ the geodesic path in $\mathfrak B$ between $u$ and $v$, composed of $k+1$ distinct vertices. Then, any path from $x$ and $y$ should at least be composed of a path from $x$ to $u \cap v_1$, a path from $v_{i-1}\cap v_i$ to $v_i\cap v_{i+1}$, and a path from $v_i\cap v_{i-1}$ to $y$. The first and last parts have length between $0$ and $\Delta_+$, while the $k-1$ others have length between $\delta_-$ (because all the $v_i$ are distinct, and so $v_{i-1}\cap v_i$ and $v_i\cap v_{i+1}$ are distinct boundary components) and $\Delta_+$.
		As a consequence, \[\delta_-\left(k-1\right) \leq d(x, y) \leq \Delta_+\left(k+1\right),\] as requested because $d(u, v) = k$.
	\end{proof}
	
	To estimate the growth rate of $\mathfrak S_\nu$, we mainly focus on the size $N_R$ of the so-called "spheres" $\mathcal S_R$. The choice of spheres rather than balls guarantees a kind of conditional independence of what is "outside" from what is "inside", which will be used to establish a geometric Markov property (Theorem \ref{th:Mark}). Nevertheless, a few times, it will be more convenient to work with $\#\mathcal B_R$. As emphasized by the following result, it amounts to the same thing.
	
	\begin{proposition}\label{prop:sphere/ball}
		Take $R\geq 0$, $P \in V^{\mathfrak B}$. Then,
		\[\#\left( \mathcal B_R(P) \setminus \mathcal S_R(P)\right) \leq \#\{\widetilde P \in \mathcal S_R(P), \mathcal B_R(P) \cap \mathcal D(\widetilde P) = \emptyset\},\] and so
		\[N_R(P) \leq \#\mathcal B_R(P)\leq 2N_R(P).\] 
	\end{proposition}
	
	\begin{proof}
		Let us prove it for $P = \rho$. Because $N_R(P) \leq \#\mathcal B_R(P)$ by definition, the second statement is a direct consequence of the first.	Consider $\mathcal B_R$ as a subtree of $\mathfrak B$. In this case, $N_R$ counts exactly the number of vertices with at most one child. The first statement means that in the binary tree $\mathcal B_R$, the number of nodes with two children is bounded by the number of leaves (\textit{i.e.} nodes without children). Build the tree by induction, adding recursively leaves to nodes of the tree, starting from $\rho$, for which the statement is true. At each step, either a child is added to a leaf, in which case the number of leaves does not change (it gains 1 and loses 1), and a new node with one child is created. Either it adds a child to a node that already has one, and simultaneously creates a leaf (the new child) and a node with two children. In both situations, the number of leaves and the number of nodes with two children increase the same, and so the statement remains true for the newly obtained tree.
	\end{proof}
	
	Finally, we state some basic properties describing the interaction between the genealogy and the spheres we defined for our study.
	\begin{proposition}\label{prop:gene}
		{Let $P, \widetilde{P} \in V^\mathfrak B$. Assume that $\widetilde{P} \in \overline{\mathcal D}(P)$. Then,
			\begin{enumerate}
				\item \label{gene:dist} $\widetilde{P} \in \mathcal  S_{d_h(\partial P^-, \partial \widetilde{P} ^-)}(P)$.
				\item \label{gene:desc} $\forall R \geq d_h(\partial P^-, \partial \widetilde{P}^-), \ \mathcal S_R(P) \cap \overline{\mathcal D}(\widetilde{P}) \neq \emptyset$.
				\item \label{gene:asc} $\forall R < d_h(\partial P^-, \partial \widetilde{P}^-), \ \mathcal S_R(P) \cap {\mathcal A}(\widetilde{P}) \neq \emptyset$.
			\end{enumerate}
		}
	\end{proposition}
	\begin{proof}
		Property \ref{gene:dist} is immediate. For \ref{gene:desc}, take a minimal element (\textit{i.e.} as deep as possible) in $\mathcal B_R(P) \cap \overline{\mathcal{D}}(\widetilde{P})$ which is non-empty as $R \geq d_h(\partial P^-, \partial \widetilde{P}^-)$. By minimality, its descendancy does not intersect $\mathcal B_R(P)$. This element is therefore in $\mathcal S_R(P)$. Property \ref{gene:asc} is obtained by a similar argument.
	\end{proof}
	
	\subsubsection{Deterministic behaviour}
	Before specifying more precisely the growth of $N_R$, let us note the existence of deterministic bounds immediately guaranteed by the geometry of pairs of pants, based on the quasi-isometry defined in Proposition \ref{prop:quasiiso}. The aim of this subsection is to establish the following lemma.
	
	\begin{lemma}[Deterministic bounds]\label{lem:majcroiss}
		For $C \geq 0$, we have \[\forall r\geq 0, \qquad 2^{\lfloor C/\Delta_+ \rfloor -3}N_r(P) \leq N_{r+C}(P) \leq e^{C+3\Delta_+} N_r(P). \tag{\ref{lem:majcroiss}.1}\label{bdet1}\] As a consequence, \[\forall R \geq 0, \qquad \frac{\ln 2}{\Delta_+}-\frac{4\ln 2}{R} \leq \frac{\ln{N_R}}{R} \leq 1+\frac{3\Delta_+}{R}.\tag{\ref{lem:majcroiss}.2}\label{bdet2}\]
		In addition, for any $P$, \[\tag{\ref{lem:majcroiss}.3}\label{bdet3} \frac{e^{-5\Delta_+}}{2} \bigop{\sum}{P' \text{ grandchild of } P}{} N_{r}(P') \leq N_{r} (P)   \leq 3\bigop{\sum}{P'\text{ grandchild of }P}{} N_r(P') .\]
	\end{lemma}
	
	The lower bound of Equation (\ref{bdet1}) is useless for the proof of Theorem \ref{th:concentration}. It is still valuable for obtaining estimates for $\alpha$, the exponential growth rate of $N_R$. In particular, it tells us that $\alpha_{\nu} > 0$, and it will be helpful in Section \ref{Surface} to obtain Theorem \ref{th:Principal3} for random lengths with support not uniformly bounded as $g \to \infty$.
	
	\begin{proof}[Proof of Lemma \ref{lem:majcroiss}]
		First, let us recall the map $\mathrm{Proj}$ introduced in Proposition \ref{prop:quasiiso}. According to this proposition, we have that \[\{P \text{ s.t. }d_{\mathfrak B}(\rho, P) \leq R/\Delta_+ - 1\} \subset \{P, d_h(\partial \rho^-, P) \leq R\} = \mathcal B_R,\] where $d_{\mathfrak B}$ is the graph distance on $\mathfrak B$, and so $\#\mathcal B_R \geq 2^{\left\lfloor {R/\Delta_+}\right\rfloor}-1\geq 2^{\left\lfloor {R/\Delta_+}\right\rfloor-1}$. The same holds for $\mathcal B_R(P)$. 
		We only prove the bound \eqref{bdet1} in the case when $P=\rho$ for simplicity of the notations. Let us start with the left side $\eqref{bdet1}$.
		Consider $P \notin \mathcal B_r$ such that its parent $\widehat P \in \mathcal B_r$ -- that is $P \in U'_r\setminus \mathcal B_r$, where $U'_r$ is defined in Definition \ref{defi:U} -- and take $\widetilde P \in \mathcal B_{C-\Delta_+}(P)$. Then, by concatenating a path from $\partial \rho^ -$ to $x \in \widehat P$ of length less than $r$, a path from $x \in \widehat P$ to some $y \in P \cap \widehat P$ of length less than $\Delta_+$ and a path from $y \in \partial P^-$ to $\widetilde P$ of length less than $C-\Delta_+$, we get that $d_h(\partial \rho^-, \widetilde P) \leq r+C$, and so $\widetilde P \in \mathcal B_{r+C}$. 
		As a consequence, \[\left(\bigop{\bigcup}{P \in U'_r, P\notin \mathcal B_r}{} \mathcal B_{C-\Delta_+}(P) \right) \subset \mathcal B_{r+C} \quad \text{and} \quad \#\left(\bigop{\bigcup}{P \in U'_r, P\notin \mathcal B_r}{} \mathcal B_{C-\Delta_+}(P) \right) \leq \#\mathcal B_{r+C} \leq 2N_{r+C},\]
		where we use Proposition \ref{prop:sphere/ball} in the last inequality.
		However, this union is disjoint: for any $\widetilde P$, there is at most one ancestor $P$ which is not in $\mathcal B_R$ but whose parent is in. In addition, any $\widehat P \in \mathcal S_r$ gives rise to a child $P \in U'_r\setminus \mathcal B_r$. Then $\#(P \in U'_r\setminus \mathcal B_r) \geq N_r$. Finally, we get that \[2^{\left\lfloor \frac{C}{\Delta_+}\right\rfloor-2}N_r \leq \bigop{\sum}{P \in U'_r, P\notin \mathcal B_r}{}  \#\mathcal B_{C-\Delta_+}(P) \leq 2N_{r+C},\] where we use our previous uniform estimate on $\# \mathcal B_{R}(P)$ with $R = C-\Delta_+$ for the first inequality.
		
		Now, let us take a look at the right-hand side of $\eqref{bdet1}$. Consider ${P} \in \mathcal S_{r + C}$, and $\gamma$ a minimal path from $\partial \rho^-$ to $P$. There exists $\widetilde{P} \in \mathcal S_{r} \cap \overline{\mathcal A}(P)$ (Property \ref{prop:gene}.\ref{gene:asc}), and $\gamma$ can be decomposed in a path $\gamma_1$ from $\partial \rho^-$ to $\partial \widetilde P^-$ and a path $\gamma_2$ from $\partial \widetilde P^-$ to ${P}$. Because $\widetilde P$ is in $\mathcal S_r$, the length of $\gamma_1$ is at least $r-\Delta_+$ by Proposition \ref{prop:elem}.\ref{elem:S}. Then, \begin{align*}
			d_h(\partial \widetilde{P}^-, \partial P^{-}) &\leq \ell(\gamma_2) = \ell(\gamma) - \ell(\gamma_1)\\ &\leq C + \Delta_+.
		\end{align*}
		In this way, $\mathcal S_{r+C} \subset \bigop{\bigcup}{P \in \mathcal S_r}{} \mathcal B_{C+\Delta_+}(P)$. 	
		Note that for any $P$, \[\#\mathcal B_{C + \Delta_+}(P) \leq e^{C + 3\Delta_+}.\] Indeed, if $\widetilde{P} \in \mathcal B_{C + \Delta_+}(P)$, then any point in $\widetilde{P}$ is at distance less than $C+2\Delta_+$ from $\partial P^-$, and at distance less than $C+3\Delta_+$ from any fixed point $x \in \partial P^-$. Because a pair of pants has area $2\pi$, we deduce that \begin{align*}2\pi \#\mathcal B_{C + \Delta_+} &\leq \mathrm{Area}(B_h(x \in \partial \rho^-, C+3\Delta_+)) \\&\leq 2\pi e^{C+3\Delta_+},
		\end{align*} since the area of a hyperbolic ball of radius $C+3\Delta_+$ in the hyperbolic plane is itself bounded by $2\pi e^{C+3\Delta_+}$.
		As a consequence,  \[N_{r+C} \leq e^{C+3\Delta_+}N_r.\]
		
		Inequality $\eqref{bdet2}$ is then an immediate consequence of $\eqref{bdet1}$ with $r=0$, $C=R$ and $N_0 = 1$.
		
		Let us now focus on Equation $\eqref{bdet3}$. Take $\widetilde{P} \in \mathcal S_r(P)$. If $\widetilde{P}$ is neither $P$ nor any of its sons, then $\widetilde{P} \in \mathcal B_r(P')$ for $P'$ a grandchild of $P$. As a consequence, \begin{align*}
			N_r(P) &\leq 3 + \bigop{\sum}{P'\text{ grandchild of }P}{} \#\mathcal B_r(P')\\ &\leq 3\bigop{\sum}{P'\text{ grandchild of }P}{} N_r(P'),
		\end{align*} where we used at the last line  Proposition \ref{prop:sphere/ball} and the fact that $N_r(P') \geq 1$ for any of the four $P'$ in the sum.
		
		Fix now $P'$ a grandchild of $P$. Any pair of pants in $\mathcal S_{r}(P')$ is also in $\mathcal B_{r+2\Delta_+}(P)$ and so  \[\bigop{\bigcup}{P'\text{ grandchild of }P}{}{\mathcal S_r(P')} \subset \mathcal B_{r+2\Delta_+}(P),\] where the union is disjoint. As a consequence, \[\bigop{\sum}{P'\text{ grandchild of }P}{}{N_r(P')} \leq \#\mathcal B_{r+2\Delta_+}(P) \leq 2e^{5\Delta_+}N_r(P),\] using both Proposition \ref{prop:sphere/ball} and Equation $\eqref{bdet1}$.
	\end{proof}
	
	Remark that similarly to the left-hand side of $\eqref{bdet1}$, we can get an upper bound \[N_{R+C}(P) \leq 2^{\lceil (C+\Delta_+)/\delta_- \rceil+1} N_R(P),\] by using the other side of the quasi-isometry inequality in Proposition \ref{prop:quasiiso} to get another upper bound on $\#\mathcal B_{C+\Delta_+}(\widetilde{P})$. Such a bound will be of no use in the proof of Theorem $\ref{th:Principal3}$, but provides information on the behaviour of $\alpha_\nu$. More about it is done in Section \ref{behavalpha}.
	
	\subsubsection{Probabilistic behaviour} Our aim is now to find more precise -- but probabilistic -- bounds on the growth rate of $N_R$. In particular, we prove Proposition \ref{prop:ssmult} and the Corollary $\ref{coro:convexp}$. 
	
	Let us begin with Proposition \ref{prop:ssmult}.
	\begin{proof}[Proof of Proposition \ref{prop:ssmult}.]
		It suffices to show that \[\forall R, r > 0,\quad N_{R+r} \leq 2\bigop{\sum}{\widetilde{P} \in U_{R}}{} N_{r+\Delta_+}(\widetilde{P}).\label{ssmult1} \tag{\ref{prop:ssmult}.2}\]
		Indeed, applying then inequality $\eqref{bdet1}$ with $C = \Delta_+$, we get the result.
		To prove inequality $\eqref{ssmult1}$, take a pair of pants $P$ in $\mathcal S_{R+r}$. Because $P \in \mathcal S_{R+r}$ then $\overline{\mathcal A}(P) \cap \mathcal S_R$ is a non-empty chain, and one can take $\widetilde{P}$ the oldest ancestor of $P$ in $\overline{\mathcal A}(P) \cap \mathcal S_R$ (it may happen that $P = \widetilde P$ when $r$ is small). Then, $\widetilde P$ has no ancestor in $\mathcal S_R$ (otherwise it would be in $\overline{\mathcal A}(P) \cap \mathcal S_R$), so it is maximal in $\mathcal S_R$. Said differently, $\widetilde P$ is in $U_R$. Consider a path from $\partial \rho^-$ to $P$ of length less than $R+r$. Because $\widetilde{P} \in \overline{\mathcal A}(P)$, $\gamma$ crosses $\widetilde P$. Then, since $\widetilde P \in \mathcal S_R$, the part of $\gamma$ between $\partial \rho^-$ and $\partial \widetilde P^-$ is of length at least $R-\Delta_+$ by Proposition \ref{prop:elem}.\ref{elem:S} and so the rest of $\gamma$ is of length at most $r+\Delta_+$. In particular, $P \in \mathcal B_{r+\Delta_+}(\widetilde{P})$. We get that $\mathcal S_{R+r} \subset \bigop{\cup}{\widetilde{P} \in U_R}{} B_{r+\Delta_+}(\widetilde{P})$ and so
		\[N_{R+r} \leq 2\bigop{\sum}{\widetilde{P} \in U_{R}}{} N_{r+\Delta_+}(\widetilde{P}),\]
		as requested.
	\end{proof}
	
	It may seem surprising to choose to state a submultiplicativity property in terms of $U_R$ instead of $\mathcal S_R$ on the right-hand side. As we will see, in the proof of Corollary \ref{coro:convexp}, it will be more convenient to use the Markov property as stated in Theorem \ref{th:Mark}. The main reason is that "everything above $U_R$" forms a complete tree -- every internal node has maximal degree -- which guarantees independence of the quantities of interest later, see Appendix \ref{geometricmarkov} for how we used it exactly.
	
	In order to obtain Corollary \ref{coro:convexp}, we need a few additional results: in what follows, we combine Property \ref{prop:ssmult} with a Markov property (Theorem \ref{th:Mark}) to get the growth rate in expectation, and establish also a supermultiplicative property to estimate the rate of convergence towards this growth rate.
	
	Let us begin with the supermultiplicative property. We use this time $U'_R$ defined in Definition \ref{defi:U} rather than looking at $U_R$. The proof is based on similar ideas to the previous ones. We first need the following technical lemma:
	\begin{proposition}\label{prop:cardinalmaj}
		Take two pairs of pants $P, \widetilde P \in V^\mathfrak B$ and $R \geq 0$. Then \[\# \left (\overline{\mathcal A}(\widetilde P) \cap \mathcal S_R(P)\right ) \leq \frac{\Delta_+}{\delta_-} + 1.\]
	\end{proposition}
	
	\begin{proof}
		Since every pair of pants in $\overline{\mathcal A}(\widetilde P) \cap \mathcal S_R(P)$ is an ancestor of $\widetilde P$, they form a finite chain for the order on the vertices of a tree. Call $Y$ the youngest one, $O$ the oldest one, and $k$ the generation gap between them (the height of $Y$ in the tree $V^\mathfrak{B}$ is the one of $O$ plus $k$). There is at most $k+1$ elements in $\overline{\mathcal A}(\widetilde P) \cap \mathcal S_R(P)$. Then, according to Property \ref{prop:elem}.\ref{elem:delta}, \[d_h(\partial P ^-, \partial O^-) + k \delta_- \leq d_h(\partial P ^-, \partial Y^-).\]
		
		Because $Y, O \in \mathcal S_R(P)$, we also have \[d_h(\partial P ^-, \partial Y^-) \leq R \leq d_h(\partial P ^-, \partial O^-) + \Delta_+,\] where the last inequality comes from Proposition \ref{prop:elem}.\ref{elem:S}, because one of the sons of $O$ is not in $\mathcal B_R$.
		This gives $k\delta_- \leq {\Delta_+}$, and so $\# \left (\overline{\mathcal A}(\widetilde P) \cap \mathcal S_r(P)\right ) \leq k+1\leq \frac{\Delta_+}{\delta_-} + 1.$ 
	\end{proof}
	
	Let us now state the supermultiplicative property.
	
	\begin{proposition}[Supermultiplicativity of the process]\label{prop:supmult}
		Fix $R \geq 0$ and recall the notation $U'_R$ from Definition \ref{defi:U}. Then
		\[\forall r \geq 0,\quad  N_{R+r} \geq \frac{e^{-4\Delta_+}}{2} \left (\frac{\Delta_+}{\delta_-} + 1 \right )^{-1} \bigop{\sum}{P \in U'_R}{} N_{r}(P). \label{supmult}\tag{\ref{prop:supmult}.1}\]
	\end{proposition}
	
	\begin{proof}
		Consider $P \in U'_R$, $\widehat P$ its parent in $\mathcal S_R$, and $\widetilde P \in \mathcal S_{r}(P)$. Then, by concatenating a path from $\partial \rho^-$ to $x \in \widehat P$ of length less than $R$, a path from $x \in\widehat P$ to some $y \in P \cap \widehat P$ of length less than $\Delta_+$ and a path from $y \in \partial P^-$ to $\widetilde P$ of length less than $r$, we get that $\widetilde P \in \mathcal B_{R+r+\Delta_+}$. 
		As a consequence, \[\bigop{\bigcup}{\widetilde P \in U'_R}{} \mathcal S_{r}(\widetilde P) \subset \mathcal B_{R+r+\Delta_+}.\]
		
		For any $P$ in $\bigop{\bigcup}{\widetilde P \in U'_R}{} \mathcal S_{r}(\widetilde P)$, we want to bound the number of pairs of pants $\widetilde P$ such that $P \in \mathcal S_{r}(\widetilde P)$. In particular, it suffices to bound $U'_R \cap \overline{\mathcal A}(P)$.	
		
		For any element in $\mathcal S_R \cap \overline{\mathcal A}(P)$, at most one of its children $\widetilde{P}$ is in $U'_R \cap \overline{\mathcal A}(P)$ by definition of $U'_R$, so $\# \left(U'_R \cap \overline{\mathcal A}(P)\right) \leq\#\left( \mathcal S_R \cap \overline{\mathcal A}(P)\right)$.
		Therefore, Lemma \ref{prop:cardinalmaj} (applied to $\widetilde P = \rho$) gives that $\# \left(U'_R \cap \overline{\mathcal A}(P)\right) \leq \left (\frac{\Delta_+}{\delta_-} + 1 \right)$ and so, 
		
		\[\bigop{\sum}{P \in U'_R}{} N_{r}(P) \leq \left(\frac{\Delta_+}{\delta_-} + 1 \right)\#\left(\bigop{\bigcup}{P \in U'_R}{} \mathcal S_{r}(P)\right) \leq  \left(\frac{\Delta_+}{\delta_-} + 1 \right)\#\mathcal B_{R+r+\Delta_+}.\]
		Then, using the fact that $\#\mathcal B_{R+r+\Delta_+} \leq 2 N_{R+r+\Delta_+}$ by Proposition \ref{prop:sphere/ball} and applying inequality $\eqref{bdet1}$ to $N_{R+r+\Delta_+}$, we get that \[\bigop{\sum}{P \in U'_R}{}  N_{r}(P) \leq 2e^{4\Delta_+}\left (\frac{\Delta_+}{\delta_-} + 1 \right ) N_{R+r},\] as requested. 
	\end{proof}

	To go further, we state a geometric Markov property for $\mathfrak S_{\nu, g}$. To do that formally, we define a notion of \textit{stopping tree}, which is a geometric equivalent of stopping times. We give a definition of the $\sigma$-algebra associated with a subtree, which roughly contains all the random variables needed to build the associated surface. In particular, we have to take into account the lengths of each boundary component of the pairs of pants in $T$, and thus all edges adjacent to $T$ (not just those in $T$).
	
	\begin{definition}\label{defi:sttr}
		Consider the following \textbf{partial order on $\rho$-rooted trees}: \[\bm{T \preceq T'} \Leftrightarrow T \subseteq T'.\]  \begin{itemize}
			\item 	For $(T, \rho)$ a subtree of $(\mathfrak B, \rho)$, we define the \textbf{$\sigma$-algebra associated to $T$} by 
			\[\bm{{\mathcal F_{T}}} = \sigma \left \lbrace \left(l_e, \tau_e\right),\ e \in E^\mathfrak B\text{ s.t. } \exists v \in V^T, \text{$e$ is adjacent to $v$ in $\mathfrak{B}$} \right \rbrace.\]
			\item We call \textbf{stopping tree} any random variable $T$ with value in the finite subtrees of $\mathfrak B$ rooted at $\rho$, such that \[\forall (T_0, \rho) \preceq (\mathcal B, \rho), \quad \mathds{1}_{T \preceq T_0} \text{is $\mathcal F_{T_0} $-measurable.}\]
			\item For $T$ a stopping tree, we set \[\bm{{\mathcal F_T}} =  \left \lbrace A\in  \mathcal F_{\mathfrak B},\ \forall (T_0, \rho) \preceq (\mathcal B, \rho), \quad A \cap \lbrace T \preceq T_0\rbrace \in \mathcal F_{T_0}\right \rbrace.\]
		\end{itemize}	
	\end{definition}
	
\begin{figure}[ht]
		\centering
		\includegraphics[width=0.5\textwidth]{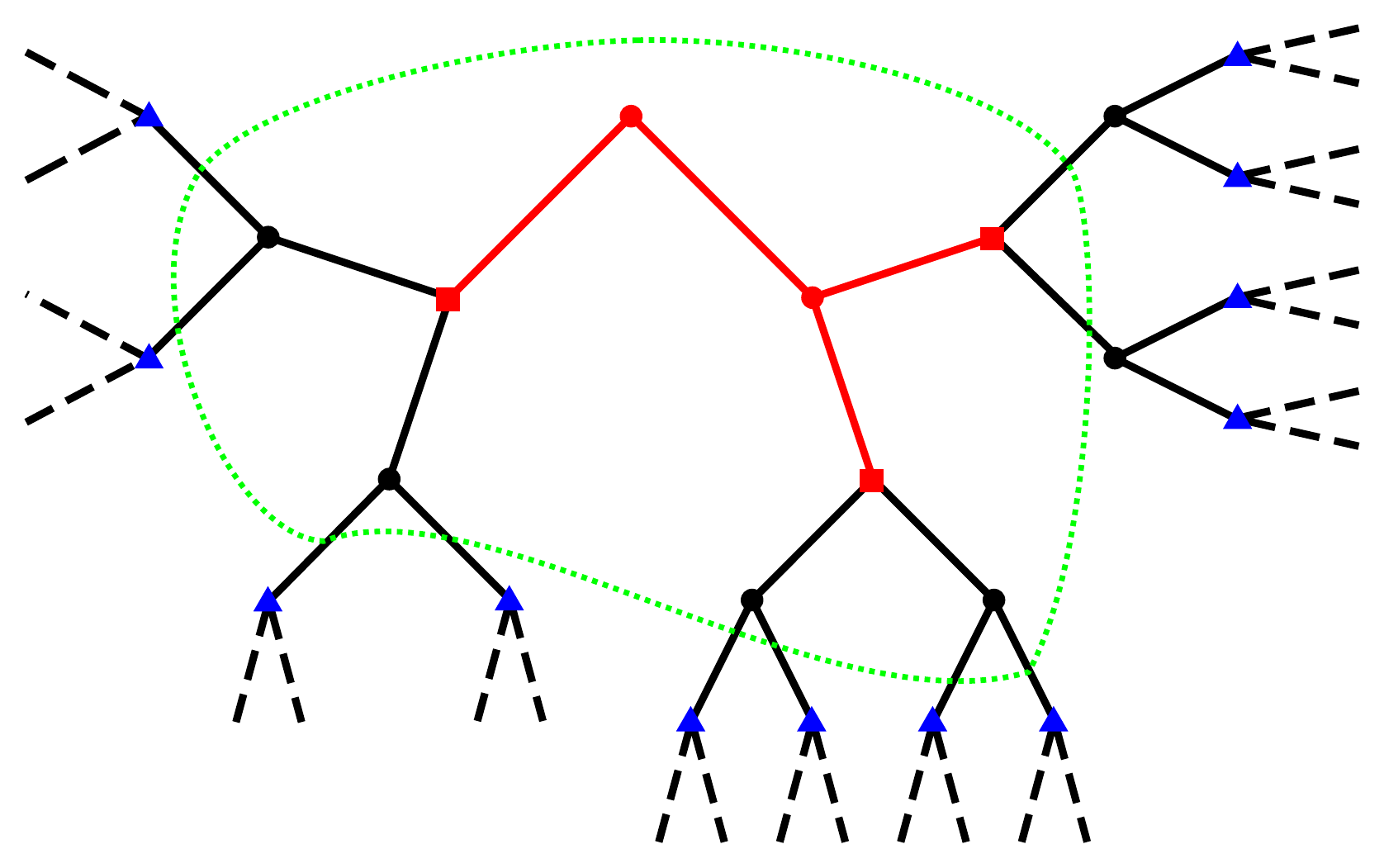}
		\caption{In the example of Figure \ref{fig:sphere}, $U_R$ corresponds to the squared vertices, and $T_R$ to the red tree. As a random variable, it depends on the weights contained in $^{\mathtt C}{}{T_R}$, surrounded in green. $\mathfrak B \setminus ^{\mathtt{C}}{}{T_R}$ is composed of 12 independent copies of $\mathfrak B$, each one being rooted in one of the blue triangular vertices.}
		\label{fig:B-T}
	\end{figure}

	To state properly the Markov Property, as in Appendix \ref{geometricmarkov}, we need to define the notion of child-tree of some tree $T$, which is just $T$ augmented with its children.
	\begin{definition}\label{defi:childtree}Let $T \preceq \mathfrak B$. 
		\begin{itemize}
			\item We define \[\bm{{^{\mathtt C}{}{T}}} = \left \lbrace P \in V^\mathfrak B,\  \text{$P$ itself or its parent is in $T$}\right\rbrace.\]
			\item Then $\bm{{\mathfrak B \setminus ^{\mathtt C}{}{T}}}$ is the subgraph induced by all the vertices which are not in $^{\mathtt C}{}{T}$ (see Figure \ref{fig:B-T}). For example, the set $U'_R$ corresponds to the leaves of $^{\mathtt{C}}{}{\mathcal B_R}$.
		\end{itemize} 
	\end{definition} 
	
	In our case, we only need the Markov property for the trees of "what is above" $U_R$ and $U'_R$:
	\begin{definition}\label{defi:treeabove} Take $R \geq 0$.
		\begin{itemize}
			\item We define \[\bm{{T_R}} = \bigop{\bigcup}{P \in U_R}{} \overline{\mathcal A}(P) = \left(\mathcal B_R \setminus \mathcal S_R \right) \cup U_R.\]
			\item We define \[\bm{ {T'_R} }= \bigop{\bigcup}{P \in U'_R}{} \overline{\mathcal A}(P) = \mathcal B_R \cup U'_R.\]
		\end{itemize}
	\end{definition}
	
	Nevertheless, it is more convenient to prove it in a more general framework, which is done in Appendix \ref{geometricmarkov} (Theorems \ref{th:Markovfaible}, \ref{th:Markovfort}). It expresses the independence of what happens inside a stopping tree from what happens outside its child-tree. We prove in Appendix \ref{geometricmarkov} that $T_R$ and $T'R$ are stopping trees (Property \ref{prop:stoptime}). Then we just focus on the following simplified version of the Markov property. 
	\begin{theorem}\label{th:Mark}
		The number $L = L(T_R)$ of leaves of $T_R$ is $\mathcal F_{T_R}$-measurable, and  $\mathfrak B \setminus ^{\mathtt C}{}{T}_R$ can be written as a disjoint union of $4L$ independent copies of $\mathcal B$, denoted $\mathfrak B_1^{T_R}, \ldots \mathfrak B_{4L}^{T_R}$ (see Figure \ref{fig:B-T}). Furthermore, for any family $(f_i)_i$ of positive functions, we have 
		\[\E\left[\bigop{\prod}{i=1}{4L(T_R)} f_i\left(\tau^{\mathfrak B_i^{T_R}}\right) \middle\vert \mathcal F_{T_R}\right] = \bigop{\prod}{i=1}{4L(T_R)} \E f_i\left(\tau^{\mathfrak B}\right).\]
		The same remains true by replacing $T_R$ by $T'_R$.
	\end{theorem}
	
	\begin{proof}
		It is a consequence of Proposition \ref{prop:stoptime}, Theorem \ref{th:Markovfort}, and of the fact that almost surely, $T_R$ (resp. $T'_R$) is finite according to the deterministic bounds of Lemma \ref{lem:majcroiss}.
	\end{proof}
	
	We now have all the framework to prove Corollary \ref{coro:convexp}. The strategy is to condition supermultiplicative and submultiplicative inequalities on what happens within distance $R$ and to use the identically distributed property of $(N_R(P), P \in V^\mathfrak B)$.

	\begin{proof}[Proof of Corollary \ref{coro:convexp}]
		First consider the submultiplicativity inequality $\eqref{ssmult}$, conditioned on $\mathcal F_{T_R}$. Using also Equation $\eqref{bdet3}$, we get that \begin{align*}\forall R, r > 0,\ \E\left[N_{R+r}\ \middle \vert \ \mathcal F_{T_R}\right] &\leq 2e^{4\Delta_+} \E\left[\bigop{\sum}{P \in V^\mathfrak B}{} \mathds{1}_{P \in U_{R}}{} N_{r}(P)\ \middle\vert\ \mathcal F_{T_R}\right]\\ &\leq 6e^{4\Delta_+} \E\left[\bigop{\sum}{P \in V^\mathfrak B}{} \mathds{1}_{\mathrm{grandparent}(P) \in U_{R}} N_{r}(P)\ \middle\vert\ \mathcal F_{T_R}\right].
		\end{align*} 
		We can now use Markov Property of Theorem \ref{th:Mark}: for any fixed $P$,  $\mathds{1}_{\mathrm{grandparent}(P) \in U_{R}}$ is $\mathcal F_{T_R}$-measurable (because so is $U_R$, the leaves of $T_R$) and $N_R(P)$ is a function of a tree in $\mathfrak B \setminus ^{\mathtt C}{}{T_R}$ : because $P$ is the grandchild of an element of $U_R$, it is precisely the root of such a tree. We get that 
		\begin{align*}
			\E\left [\bigop{\sum}{P \in V^\mathfrak B}{} \mathds{1}_{\mathrm{grandparent}(P) \in U_{R}}{} N_{r}(P)\ \middle \vert \ \mathcal F_{T_R}\right] &= \bigop{\sum}{P \in V^\mathfrak B}{} \mathds{1}_{\mathrm{grandparent}(P) \in U_{R}}{} \E N_{r}(P)\\ &\leq 4\#U_R \E N_{r} \leq 4N_R \E N_{r},
		\end{align*}
		Combining the two previous inequalities and taking expectations on both sides, we obtain that
		\[\forall R, r \geq 0,\quad\E N_{R+r} \leq 24e^{4\Delta_+} \E N_R \E N_r.\]
		In this way, the map $R \in \Rp \mapsto  \ln\E N_R +  4\Delta_+ +\ln 24$ is subadditive. As a consequence of Fekete's Lemma, \[\alpha_\nu = \inf_{R \geq 0} \frac{\ln \E N_R+4{\Delta_+}+\ln 24}{R} = \lim_{R \to \infty} \frac{\ln \E N_R + 4{\Delta_+} + \ln 24}{R} = \lim_{R \to \infty} \Sigma_R,\]
		where we recall about the notation $\Sigma_R = \frac{\ln\E N_R}{R}$.
		In particular, \[\alpha_\nu - \Sigma_R \leq \frac{4{\Delta_+} + \ln 24}{R}.\]
		
		Applying the same reasoning to the supermultiplicative inequality $\eqref{supmult}$ we get that \begin{align*}
			\forall R, r \geq 0,\quad \E\left[N_{R+r}\ \middle|\ \mathcal F_{T'_R}\right] &\geq \frac{e^{-4\Delta_+}}{2} \left(\frac{\Delta_+}{\delta_-}+1\right )^{-1} \E\left[\bigop{\sum}{P \in U'_R}{}N_{r}(P)\ \middle|\ \mathcal F_{T'_R}\right] \\&\geq \frac{e^{-9\Delta_+}}{4} \left (\frac{\Delta_+}{\delta_-}+1\right )^{-1} \bigop{\sum}{P \in V^\mathfrak B}{} \mathds{1}_{\mathrm{grandparent}(P) \in U'_{R}} \E\left[N_{r}(P)\ \middle|\ \mathcal F_{T'_R}\right].\end{align*}
		
		We get as previously that the function $R \in \Rp \mapsto  \ln\E N_R -  \ln\left (\frac{\Delta_+}{\delta_-}+1 \right ) - 9\Delta_+-\ln 4$ is superadditive, and consequently, \[\sup_{R \geq 0} \frac{\ln\E N_R - \ln\left (\frac{\Delta_+}{\delta_-}+1 \right )- 9\Delta_+-\ln 4}{R} = \lim_{R \to \infty} \Sigma_R = \alpha_\nu.\]
		
		In particular, \[-\frac{\ln \left( \frac{\Delta_+}{\delta_-}+1 \right)+ 9\Delta_++\ln 4}{R} \leq \alpha_\nu - \Sigma_R \leq \frac{4\Delta_+ + \ln 24}{R}.\]
		
		All together, \begin{align*} \left\lvert \alpha_\nu - \Sigma_R \right\rvert &\leq \frac{\max \left(4\Delta_+ + \ln 24, \ln \left( \frac{\Delta_+}{\delta_-}+1 \right)+ 9\Delta_++\ln 4\right)}{R}\\&\leq \frac{15\Delta_+}{R},\end{align*}
		because $\ln 4 \leq 2\Delta_+$ and $\ln\left(\frac{\Delta_+}{\delta_-}+1\right) \leq 4\Delta_+$ according to Proposition \ref{prop:Ddelta}.
	\end{proof}
	
	In the following, we simply write $\alpha$ for $\alpha_\nu$.
	
	\subsection{Almost sure convergence}\label{convas}
	
	\subsubsection{Proof of Theorem \ref{th:concentration}}
	In this section, we prove Theorem \ref{th:concentration}, starting from the convergence of $\Sigma_R$ around $\alpha$. The proof is inspired by what is done in \cite[Section 3.2]{Bordenave}. The idea is that if we have a concentration of $\frac{\ln N_R}{R}$ around its expected value, and we know that $\E N_R^2$ does not grow too fast, we can prove that $\Sigma_R$ and $\E \frac{\ln N_R}{R}$ are close enough to give us a concentration of $\frac{\ln N_R}{R}$ around $\Sigma_R$, and so around $\alpha$. To do so, we first state the two following technical lemmas.
	
	\begin{lemma}\label{lem:Majca} We have that
		\[\forall \eta > 0, \exists R_1, \forall R \geq R_1, \quad \E N_R^2\leq \left(1 + \eta\right)^R \E[N_R]^2.\]
		We can take $R_1 = \left(\frac{153\Delta_+}{\ln\left(1+\eta\right)}\right)^2$.
	\end{lemma}
	
	\begin{proof}
		Fix $\eta > 0$. We want to show that for $R$ large enough, $\E N_R^2 \leq (1+\eta)^R\ \E[N_R]^2$, that is \[\gamma_R = \frac{\ln{\E N_R^2}}{R} - 2 \Sigma_R \leq \ln(1+\eta).\]
		
		Consider the square of the submultiplicativity inequality $\eqref{ssmult}$, for $R = r$. This gives \[\forall R > 0,\quad  N_{2R}^2 \leq 4e^{8\Delta_+} \left( \bigop{\sum}{P \in U_{R}}{} N_{R}^2(P) + \bigop{\sum}{P \neq \widetilde{P} \in U_{R}}{} N_{R}(P)N_{R}(\widetilde{P}) \right).\]
		
		Combined with inequality $\eqref{bdet3}$, we obtain that
		\begin{multline*}N_{2R}^2 \leq 36e^{8\Delta_+} \bigop{\sum}{P \in U_{R}}{} \left(\bigop{\sum}{P'\text{ grandchild of }P}{} N_r(P')\right)^2\\ + 36e^{8\Delta_+} \bigop{\sum}{P \neq \widetilde{P} \in U_{R}}{}  \left(\bigop{\sum}{P'\text{ grandchild of }P}{} N_r(P')\right)\left( \bigop{\sum}{P'\text{ grandchild of }\widetilde P}{} N_r(P')\right).\end{multline*}
		
		Let us first focus on the sum on the first line.
		\begin{align*}
			\bigop{\sum}{P \in U_{R}}{} \left(\bigop{\sum}{P'\text{ grandchild of }P}{} N_r(P')\right)^2 &= \bigop{\sum}{P \in V^\mathfrak B}{} \mathds{1}_{\mathrm{grandparent}(P) \in U_{R}}{}  N_{r}^2(P) \\&+ 2\bigop{\sum}{\substack{P \neq P' \in V^\mathfrak B\\\text{with the same grandparent}}}{} \mathds{1}_{\mathrm{grandparent}(P) \in U_{R}}{} N_{r}(P)N_{r}(P').
		\end{align*}
		As in the proof of Corollary \ref{coro:convexp}, for fixed $P \neq P'$ with the same grandparent,  $\mathds{1}_{\mathrm{grandparent}(P) \in U_{R}}$ is $\mathcal F_{T_R}$-measurable and $N_R(P)$ (respectively $N_R(P)N_R(P')$) is a function of a tree (respectively the product of two functions of distinct trees) in $\mathfrak B \setminus ^{\mathtt C}{}{T_R}$. We therefore can use the Markov inequality \ref{th:Mark} to get
		\begin{align*}&\E\left[\bigop{\sum}{P \in U_{R}}{} \left(\bigop{\sum}{P'\text{ grandchild of }P}{} N_r(P')\right)^2\middle|\ \mathcal F_{T_R}\right]  \\&= \bigop{\sum}{P \in V^\mathfrak B}{} \mathds{1}_{\mathrm{grandparent}(P) \in U_{R}}{}  \E N_{r}^2(P) + 2\bigop{\sum}{\substack{P \neq P' \in V^\mathfrak B\\\text{with same grandparent}}}{} \mathds{1}_{\mathrm{grandparent}(P) \in U_{R}}{} \E N_{r}(P)\E N_{r}(P').
		\end{align*}
		In particular,
		\[\E\left[\bigop{\sum}{P \in U_{R}}{} \left(\bigop{\sum}{P'\text{ grandchild of }P}{} N_r(P')\right)^2\right]  \leq 4\#U_{R}  \E N_{r}^2 + 2 \binom{4}{2}\#U_R \E\left[N_{r}\right]^2.\]
		Similarly, for the second line of the previous upper bound, we get that \begin{align*}&\E\left[\bigop{\sum}{P \neq \widetilde{P} \in U_{R}}{}  \left(\bigop{\sum}{P'\text{ grandchild of }P}{} N_r(P')\right)\left( \bigop{\sum}{P'\text{ grandchild of }\widetilde P}{} N_r(P')\right) \right]  \\&\leq \#\left\{(P_1, P_2), \text{with $\neq$ grandparents in $U_R$}\right\} \E\left[N_{r}\right]^2\\&\leq 16\left(\#U_R^2 - \#U_R\right)\E[N_R]^2.
		\end{align*}
		Finally, we get that
		\[\E N_{2R}^2 \leq 36e^{8\Delta_+} \left(4\E \#U_{R} \E N_R^2 +  16 \E\# U_R^2\E[N_R]^2\right) \leq 720e^{8\Delta_+} \E N_R^2\E[N_R]^2,\] where we use that $\#U_R \leq N_R \leq N_R^2$.
		Take the logarithm and divide by $2R$ to get:
		\begin{align*}
			\gamma_{2R} + 2\Sigma_{2R} &\leq \frac{7\Delta_+}{R} + \frac{\gamma_R}{2} + 2\Sigma_R. \end{align*}
		because $\frac{\ln 720}{2} \leq 3\Delta_+$.	
		Using the rate of convergence of $\Sigma_R$ towards $\alpha$ given by inequality $\eqref{ratealpha}$, this leads to \begin{align*}
			\gamma_{2R} &\leq \frac{\gamma_R}{2} + \frac{7\Delta_+}{R} + 2(\Sigma_R - \alpha) + 2(\alpha-\Sigma_{2R})\\&\leq \frac{\gamma_R}{2} + \frac{52\Delta_+}{R}.
		\end{align*}
		Iterating this, we obtain \[\forall k \in \N, \quad \gamma_{2^kR} \leq \frac{\gamma_{2^{k-1}R}}{2} + \frac{52\Delta_+}{2^{k-1}R}\leq \frac{\gamma_R}{2^k} + \frac{52k\Delta_+}{2^{k-1}R}.\]
		Take $0<\eta<1$. Then, thanks to the deterministic bounds $\eqref{bdet1}$, we know that for $R \geq 1$, $
		\gamma_R \leq \frac{\ln \E N_R^2}{R} \leq 8 \Delta_+$.
		As a consequence, \[\forall k \geq 2, \forall R\geq 2^k, \quad \gamma_{R} \leq \frac{8\Delta_+}{2^{k}} + \frac{52k\Delta_+}{2^{k}} \leq \frac{108\Delta_+}{2^{k/2}}\] because $k2^{-k} \leq 2^{1-k/2}$ and $2^{k} \geq 2^{k/2+1}$.
		Take $R \geq R_1 = \left (\frac{153\Delta_+}{\ln\left(1+\eta\right)}\right)^2 > 2^{\left\lceil 2\log_2\left(\frac{108\Delta_+}{\ln\left(1+\eta\right)}\right)\right\rceil}$. Then, we have shown that \[ \gamma_R \leq \ln(1+\eta),\] as requested.
	\end{proof}
	Note that instead of bounding $k2^{-k}$ by  $2^{1-k/2}$, we could have worked harder to obtain a bound of the form $C2^{-k(1-\varepsilon)}$ for any $\varepsilon>0$. All in all, it would have led to better fluctuations on the growth rate of the surfaces (roughly of order $\ln^{1/2} g$).
	
	\begin{lemma}[{Concentration of $\ln N_R$}]\label{lem:cMDiarmid}
		There exists some numerical constant $C$, independent of all the parameters of the model, such that  \[\forall \varepsilon > 0, \forall R \geq 2\delta_-, \quad \p\left( \left \lvert \frac{\ln N_R}{R} - \E \frac{\ln N_R}{R}\right\rvert \geq \varepsilon \right) \leq 2\exp\left(-e^{-C\Delta_+}\varepsilon^2R \right).\]
	\end{lemma}
	
	\begin{proof}[Proof of Lemma \ref{lem:cMDiarmid}]
		By Proposition \ref{prop:quasiiso}, $N_R$ counts pairs of pants which are at most at depth $\frac{R}{\delta_-}$ so is determined by weights on edge between pairs of pants of height less than $n = \left \lceil \frac{R}{\delta_-}\right \rceil$. We want to apply McDiarmid's inequality (recalled in Proposition \ref{prop:McD}) to this function, and to the random variables $(X_i)_{0\leq i \leq n}$, where $X_i$ is the $2^i$-uplet of twists parameters associated with edges between pairs of pants of height $i$ and of height $i+1$.
		
		Let us fix $i$ between $0$ and $n$, and $X = (X_0, \ldots X_i, \ldots, X_n)$, $X' = (X_0, \ldots X'_i, \ldots, X_n)$, with $X'_i$ an independent copy of $X_i$. Let $\gamma$ be a geodesic path from $\partial \rho^-$ to some fixed pair of pants $P$ in $\mathfrak S_\nu(X)$. It may happen that $\gamma$ (or its projection on $\mathfrak B$) leaves several times the geodesic of the tree $\mathfrak B$ between $\rho$ and $P$ to visit pairs of pants at height $i$. Let $P'$ be a pair of pants where such divergence happens (Note that $\gamma$ may also diverge without visiting pairs of pants at height $i$, but this is not relevant for us). Then, we modify $\gamma$ so that as soon as it enters $P'$, it directly goes (with the shortest path inside $P'$) to the last point of $P'$ it visits, without exiting $P'$ in the meantime. Each of these modifications adds a factor at most $\Delta_+$ to the length of the path. Let $\gamma'$ be the path we finally obtain, and $K_i$ the number of modifications that have been made. Now, $\gamma'$ is of length at most $\ell(\gamma) + K_i\Delta_+$, and crosses at most two pairs of pants at height $i$ or $i+1$ (if $P$ is deep enough, the two pairs of pants between $\rho$ and $P$ in the tree).
		Let us now consider the same path $\gamma'$ but in the surface $\mathfrak S_\nu(X')$. The length of the part of $\gamma'$ before entering the unique pair of pants at height $i$ is not modified, and similarly after the last time it leaves height $i+1$. By joining these two parts with a path of length at most $2\Delta_+$ (to cross the two pairs of pants at height $i$ and $i+1$), we get a path from $\rho$ to $P$ in $\mathfrak S_\nu(X')$ of length at most $\ell(\gamma) + (K_i+1)\Delta_+ \leq R + (K_i+1)\Delta_+$.
		
		We now claim that $K_i \leq K_0 = 92(1/\delta_-+\Delta_+)^2$ for any $P$. If so, we have that \[\mathcal S_R(X) \subset \mathcal B_R(X) \subset \mathcal B_{R+K_0\Delta}(X').\]
		As a consequence, \[N_R(X) \leq \#B_{R+K_0\Delta_+}(X') \leq 2N_{R+K_0\Delta_+}(X') \leq 2e^{(K_0+3)\Delta_+}N_r(X'), \]
		using the Proposition \ref{prop:sphere/ball} and the deterministic bound $\eqref{bdet1}$ for the two last inequalities.
		At the end, we get that \[\lvert \ln N_R(X) - \ln N_R(X')\rvert \leq c_i =  (K_0+4)\Delta_+.\] 
		McDiarmid's inequality finally tells us that \[\p\left(\lvert \ln N_R - \E\ln N_R\rvert \geq \varepsilon\right) \leq 2\exp\left(-\frac{2\varepsilon^2}{\left(\left\lceil\frac{ R}{\delta_-}\right\rceil+1\right) (K_0+4)^2\Delta_+^2}\right),\] and therefore
		\begin{align*}\p\left(\left\lvert \frac{\ln N_R}{R} - \E \frac{\ln N_R}{R}\right\rvert \geq \varepsilon\right) &\leq 2\exp\left(-\frac{2\varepsilon^2R^2}{\left(\left\lceil\frac{ R}{\delta_-}\right\rceil+1\right) (K_0+4)^2\Delta_+^2}\right)\\&\leq 2\exp\left(-\frac{\delta_-\varepsilon^2R}{\left(95(\delta_-^{-1}+\Delta_+)^2\Delta_+\right)^2} \right), \end{align*} when $R>\delta_-$, since $4<3(\delta_-^{-1}+\Delta_+)^2$. Remark that according to Proposition \ref{prop:Ddelta}, $\delta_- \geq e^{-\Delta_+}$, and so because $\Delta_+$ is uniformly bounded from below, $\frac{\delta_-}{(95(\delta_-^{-1}+\Delta_+)^2\Delta_+)^2} \geq e^{-C\Delta_+}$ for some numerical $C$ large enough, independent of any other variables. This is enough provided that we indeed have $K_i \leq K_0$ for any $P$. If $i\leq K_0$, this is immediate because $K_i\leq i$. Suppose now that $i>K_0$.
		According to Proposition \ref{prop:quasiiso}, $\mathfrak S_\nu$ is $\left(\Delta_+ + \delta_-^{-1}, 1\right)$-quasi-isometric to $\mathfrak B$ which is 0-hyperbolic. As a consequence, $\gamma$ projects into what is called a $\left(\Delta_+ + \delta_-^{-1}, 1\right)$-quasi-geodesic in hyperbolic geometry. According to \cite{Gouezel}, this quasi-geodesic must be in the $92\left(\Delta_+ + \delta_-^{-1}\right)^2$-neighbourhood of the true geodesic in the tree. As a consequence, to reach a pair of pants at height $i$, $\gamma$ can only diverge at pairs of pants at height between $i-K_0$ and $i$, so it happens at most $K_0$ times.
	\end{proof}
	
	Note that the estimate of $c_i$ does not seem to be good for $i$ large. In fact, when we say that $N_R$ is affected by all twists up to height $\frac{R}{\delta_-}$, we are considering the case where the path between the root and the pairs of pants at height $h \approx \frac{R}{\delta_-}$ is roughly as short as possible. However, in general, few of the $2^h$ pairs of pants at height $h$ are in this case. As a consequence, changing all the twists at the height of about $\frac{R}{\delta_-}$ will have a real effect in very few branches of $\mathcal S_R$, which is not taken into account in the estimation of $c_i$ we have done here. We can imagine that a better estimation could improve the concentration inequality. If we were also able to get a better upper bound on the uniform $K_0$ (for example with more information on $\nu$), we could hope to get rid of the terms in $\delta_-$ in the proof. As we discussed in the introduction, this would allow us to improve our result on the diameter with fewer restrictions on the support of the length law.
	
	Let us now see how it is sufficient to conclude to Theorem \ref{th:concentration}.
	
	\begin{proof}[Proof of Theorem \ref{th:concentration}] 
		According to Paley-Zigmund inequality (recalled in Proposition \ref{prop:PZ}), we have  \[\p\left(\frac{\ln N_R}{R} \geq \Sigma_R - \frac{1}{R}\right) = \p\left(N_R \geq e^{-1}\E N_R \right) \geq \left(1-e^{-1}\right)^2 \frac{\E[N_R]^2}{\E N_R^2}.\] Let us write $\frac{\E[N_R]^2}{\E N_R^2} = \left(1+\eta_R\right)^{-R}$ for some $\eta_R \geq 0$.
		Then, for $0< \varepsilon <1$, for $R \geq \delta_-$, \begin{align*}
			\p\left(\Sigma_R-\frac{1}{R} \leq \frac{\ln N_R}{R} < \E \frac{\ln N_R}{R}+\varepsilon \right) &\geq \p\left (\Sigma_R-\frac{1}{R} \leq \frac{\ln N_R}{R}, \left|\frac{\ln N_R}{R}-\E \frac{\ln N_R}{R}\right| < \varepsilon\right) \\&\geq \p\left(\Sigma_R-\frac{1}{R} \leq \frac{\ln N_R}{R}\right) - \p\left(\left\lvert \frac{\ln N_R}{R}-\E \frac{\ln N_R}{R} \right\rvert \geq \varepsilon\right) \\&\geq \left(1-e^{-1}\right)^2\left(1+\eta_R\right)^{-R} - 2\exp\left(-e^{-C\Delta_+}\varepsilon^2R\right)
		\end{align*}
		using Lemma \ref{lem:cMDiarmid}.
		
		This is positive as soon as \[\ln(1+\eta_R) \leq e^{-C\Delta_+} \varepsilon^2 - \frac{2}{R} < e^{-C\Delta_+} \varepsilon^2 - \frac{2\ln\left(\frac{\sqrt{2}e}{e-1}\right)}{R}.\] In particular, if $R \geq \frac{4e^{C\Delta_+}}{\varepsilon^2}$, it is true whenever $\ln(1+\eta_R) \leq \frac{e^{-C\Delta_+}\varepsilon^2}{2}$.
		Thus, from Lemma \ref{lem:Majca}, we deduce that this is positive as soon as \[R \geq \left(\frac{306\Delta_+}{e^{-C\Delta_+}\varepsilon^2}\right)^2 > \frac{4e^{C\Delta_+} }{\varepsilon^2}.\] It is enough to take $R \geq {R}_0 = \frac{e^{c\Delta_+}}{\varepsilon^4}$ for some well-chosen numerical constant $c>0$, not depending on any other variable.
		So with positive probability,  $\Sigma_R - 1/R \leq \frac{\ln N_R}{R} \leq \E \frac{\ln N_R}{R} + \varepsilon$. Since $\E \frac{\ln N_R}{R}$ and $\Sigma_R - 1/R$ are deterministic quantities, we deduce that, for $R \geq {R}_0 \geq 1/\varepsilon$, \[\Sigma_R - 2\varepsilon \leq \E \frac{\ln N_R}{R} \leq \Sigma_R,\] the second inequality coming from concavity of logarithm.
		Thus, for $R \geq {R}_0$,
		\begin{align*}\p\left(\left\lvert \frac{\ln N_R}{R} - \Sigma_R \right\rvert \geq 3\varepsilon\right) &\leq \p\left(\left\lvert \frac{\ln N_R}{R} - \E \frac{\ln N_R}{R} \right\rvert \geq \varepsilon\right) \\&\leq 2\exp\left(-{e^{-C\Delta_+}}\varepsilon^2R \right),\end{align*} according to Lemma \ref{lem:cMDiarmid}.
		
		Furthermore, Corollary \ref{coro:convexp} says that for $R \geq \frac{15\Delta_+}{\varepsilon}$, $\lvert \Sigma_R - \alpha \rvert \leq \varepsilon$ so finally, \[\forall R \geq R_0 = \frac{e^{c\Delta_+}}{\varepsilon^4}, \quad	
		\p\left(\left\lvert \frac{\ln N_R}{R} - \alpha \right\rvert \geq 4\varepsilon\right) \leq 2\exp\left(-e^{-C\Delta_+} \varepsilon^2R \right),\]
		which gives the result for another value of $C$.
	\end{proof}
	
	\subsubsection{What about \texorpdfstring{$\widehat{\mathfrak S}_\nu$}{Ŝ\_ν} ?}
	To conclude this part, we give a corollary which describes what happens for $\widehat{\mathfrak S}_\nu$, the surface built with the whole tree $\mathfrak T_3$. We write $\widehat{\mathcal B}_R, \widehat{\mathcal S}_R, \widehat{N}_R$ for the equivalents on $\widehat{\mathfrak S}_\nu$ of $\mathcal B_R, \mathcal S_R, N_R$
	\begin{corollary}[{Behaviour for $\widehat{\mathfrak S}_\nu$}]\label{coro:surfacecomplete}
		For all $0<\varepsilon <1$, and $R \geq R_0$,
		\begin{align*}
			\p\left(\widehat{N}_R \leq e^{(\alpha-\varepsilon) R}\right)&\leq 2 \exp\left(-e^{-C\Delta_+}\varepsilon^2\left(R+\frac{\ln 2}{\alpha-\varepsilon}\right)\right),\\ \p\left(\widehat{N}_R \geq e^{(\alpha+\varepsilon) R}\right) &\leq 6 \exp\left(-e^{-C\Delta_+}\varepsilon^2\left(R-\frac{\ln(3/2)}{\alpha+\varepsilon}\right)\right).
		\end{align*}
	\end{corollary} 
	
\begin{figure}[ht]
		\centering
		\includegraphics[width=\textwidth]{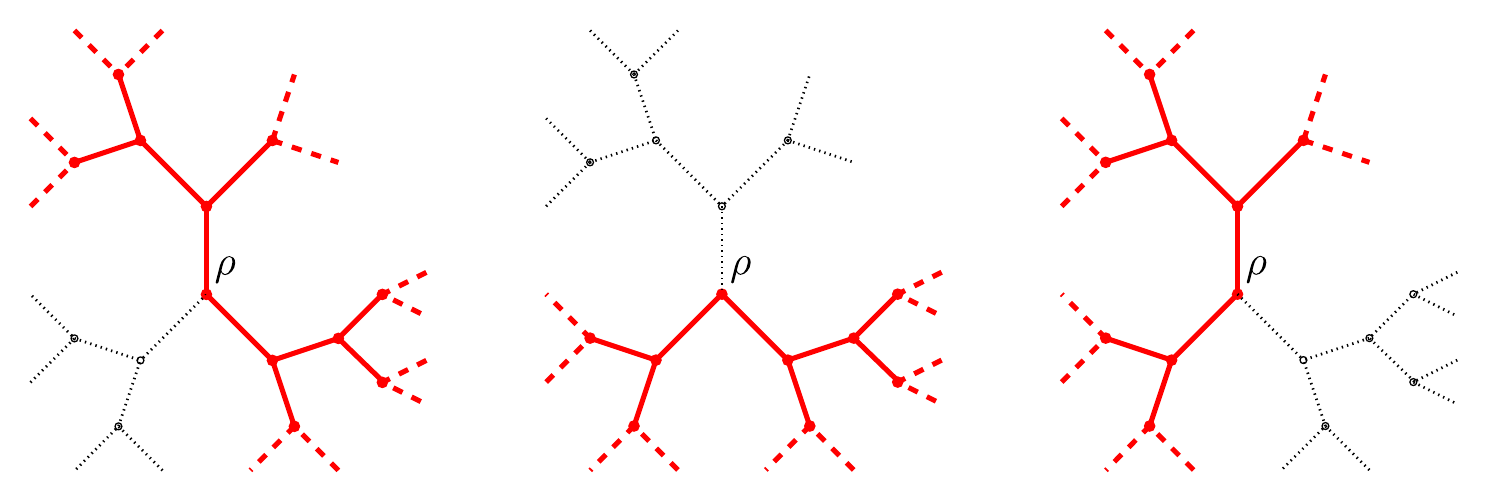}
		\caption{To estimate the number of pairs of pants within distance $R$ from $\rho$ in $\widehat{\mathfrak S}_\nu$, we use our estimates in each of the bold copies of $\mathfrak S_\nu$, giving 3 copies of $N_R$. Indeed, $2\widehat{N}_R$ is equal to the sum of these 3 (non-independent) copies.}
		\label{fig:treelike_sum}
	\end{figure}

	\begin{proof}
		We can see $2\widehat{N}_R$ as the sum of 3 (non-independent) copies of $N_R$ (see Figure \ref{fig:treelike_sum}). With this and Theorem \ref{th:concentration} in mind, we get
		\begin{align*}
			\p\left(\widehat{N}_R \leq e^{(\alpha-\varepsilon) R}\right)&= \p\left(2\widehat{N}_R \leq 2e^{(\alpha-\varepsilon) R}\right)\\&\leq \p\left({N}_R \leq 2e^{(\alpha-\varepsilon) R}\right) \\&\leq 2 \exp\left(-e^{-C\Delta_+}\varepsilon^2\left(R+\frac{\ln 2}{\alpha-\varepsilon}\right)\right), \end{align*} and
		\begin{align*} \p\left(\widehat{N}_R \geq e^{(\alpha+\varepsilon) R}\right) &= \p\left(2\widehat{N}_R \geq 2e^{(\alpha+\varepsilon) R}\right)\\ &\leq 3\p\left({N}_R \geq \frac{2}{3}e^{(\alpha+\varepsilon) R}\right)\\&\leq 6 \exp\left(-e^{-C\Delta_+}\varepsilon^2\left(R-\frac{\ln(3/2)}{\alpha+\varepsilon}\right)\right),\end{align*}
		as requested.
	\end{proof}
	
	\subsection{Asymptotic behaviour of \texorpdfstring{$\alpha$}{α}}\label{behavalpha}
	In this section, we discuss the asymptotic behaviour of $\alpha$ in the particular case where all the boundary lengths are equal to a deterministic value $2l>0$ -- \textit{i.e.} when $\nu_l$ is a Dirac mass in $2l$. In this case, and when $\nu_t$ is the twists distribution, we write $\bm{{\alpha_{l\otimes\nu_t}}}$ for $\alpha_\nu$ by abuse of notation.
	
	Note that Lemma \ref{lem:majcroiss} and the remark following it directly tell us that \[\frac{\ln 2}{\Delta_+} \leq \alpha_\nu\leq \min\left(1, \frac{\ln 2}{\delta_-}\right),\]
	and so \[\alpha_{l\otimes\nu_t} \underset{l \to 0}{\longrightarrow} 0,\] independently of the law of twists $\nu_t$. This is not surprising: when $l \to 0$, it is well known that $\delta_- \to \infty$ (see the expression of $\delta_-$ or the collar lemma in Appendix \ref{collar} for example), which amounts to saying that it takes "a long time before discovering a new pair of pants", hence the small growth rate.
	
	Nevertheless, these inequalities do not say much when $l \to \infty$. In the special case where $\nu_t$ is the Dirac mass in $0$, \textit{i.e.} twists almost surely null, $\alpha$ can be interpreted as a nice geometric quantity: as stated in \cite[Theorem 2]{BCP_mindiam}, it is the critical exponent of the Fuchsian group generated by the reflections in three non-consecutive sides of a right angled hyperbolic hexagon such that the three non-consecutive sides all have length $l$. This interpretation allows the authors to show that \[\lim_{l \to \infty} \alpha_{l \otimes 0} = 1.\]
	
	The aim of this section is to extend this last result to the case of uniformly chosen twists, that is, to obtain Theorem \ref{th:balpha}.
	More precisely, we already know that $\limsup_{l \to \infty }\alpha_{l\otimes\nu_t} \leq 1$ by inequality \eqref{bdet2}. Let us therefore show that the volume -- and consequently $N_R$ -- grows fast enough. To do so, we define a set of "good" pairs of pants, which will roughly be the ones from which we can easily obtain a good growth.
	
	\subsubsection{Notations}
	
	Consider $P \in V^\mathfrak B$, whose boundary components are $\partial P^-$, $\partial P^+_1$ and $\partial P^+_2$, all of length $2l$. To work properly, we first explicit the arc-length parameterizations $\partial P^-(t)$, $\partial P^+_1(t)$ and $\partial P^+_2(t)$ of the boundary components of $P$.
	\begin{figure}
		\centering
		\subfloat[0.4\textwidth][Parametrization of the boundary.]
		{
			\centering
			\includegraphics[width=0.35\textwidth]{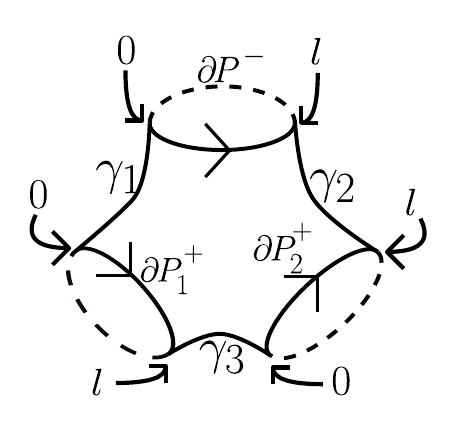}
			\label{fig:boundparam}
		}
		\hfill
		\subfloat[0.4\textwidth][The curve $\widetilde{\gamma}$ cut the pair of pants into two.]
		{
			\centering
			\includegraphics[width=0.35\textwidth]{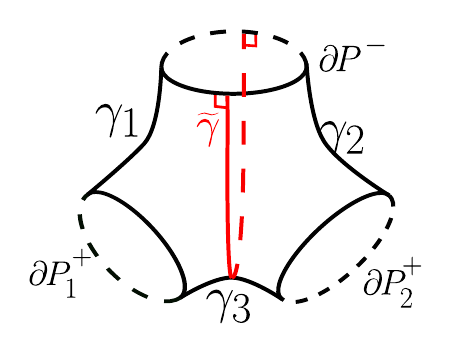}
			\label{fig:halfpant}
		}
		\caption{A pair of pants.}
		\label{fig:parampants}
	\end{figure}
	
	\begin{definition}\label{defi:curves}
		We proceed as follows: 
		\begin{itemize}
			\item There is a unique geodesic path $\bm{{\gamma_1}}$ between $\partial P^-$ and $\partial P_1^+$ which is orthogonal to both of them. We choose ${\partial P^-(0)}$ to be the point of $\partial P^-$ on this geodesic, and similarly, ${\partial P_1^+(0)}$ is the point of $\partial P_1^+$ on this geodesic.
			\item On the same principle, there is a unique geodesic path $\bm{{\gamma_2}}$ between $\partial P^-$ and $\partial P_2^+$ orthogonal to both of them, which joins $\partial P^-(l)$ and $\partial P_2^+(l)$.
			\item Consequently, the unique orthogonal geodesic path $\bm{{\gamma_3}}$ between $\partial P_1^+$ and $\partial P_2^+$ joins $\partial P_1^+(l)$ and $\partial P_2^+(0)$.
			\item About \textbf{the orientation} of the parameterization, we require that $\partial P^-\left([0, l]\right)$, $\partial P_1^+\left([0, l]\right)$, $\partial P_2^+\left([0, l]\right)$ belongs to the same hyperbolic hexagon, and that when two pairs of pants $P$ and $\widetilde{P}$ are glued along $\partial P_i^+$ and $\partial \widetilde{P}^-$ with twist $\tau$, we have \[\partial P_i^+(t)=\partial \widetilde{P}^-(t-\tau).\] 
			Be careful that this is not the usual convention (and in particular not what is done in \cite[Chapter 3]{Buser}): one way to do so is to take $\partial P^-$ directly oriented, and the $\partial P_i^+$ indirectly oriented. 
		\end{itemize}
	\end{definition}
	As soon as the numbering does not depend on the randomness of the weights, changing it does not change the distribution of the resulting surface. This can be extended even when the lengths are random, showing that the notion of twist is well-defined.
	
	With this setting, we can define the notion of "good" pairs of pants that will ensure a good growth rate later.
	
	\begin{definition}
		Let $\bm{{\mathrm{Good}_R}}$ be the set of pairs of pants in $\mathcal S_R$ that are leaves of $\mathcal B_R$ or that are "intersected" by the ball $\mathcal B_R$ at a point in $\partial P^-\left(\left[\frac{l}{2}-1; \frac{l}{2}+ 1\right]\right)$: \[\mathrm{Good}_R = \left \lbrace P\in \mathcal S_R,\ \mathcal D(P)\cap\mathcal B_R = \emptyset \text{ or }  \partial P^-\left(\left[\frac{l}{2}-1; \frac{l}{2}+ 1\right]\right) \cap B_h(\partial \rho^-, R)\neq \emptyset\right\rbrace,\] where we recall that $B_h(\partial \rho^-, R)$ is the hyperbolic ball of radius $R$ centered at $\partial \rho^-$.
	\end{definition}
	
	\subsubsection{Proof of Theorem \ref{th:balpha}}
	
	The main idea is to proceed by induction on $k \geq 1$ to bound from below $N_{kl}$, or equivalently $\#\mathcal B_{kl}$. To do so, we first prove that the injectivity radius of $\widehat{\mathfrak S}_l$ is exactly $l$, and then use $\Good_{kl}$ to find enough distinct balls of radius $l$ (whose volumes are the same as in the hyperbolic plane) to estimate $\#\mathcal B_{(k+1)l}$. 
	\begin{proposition}[Injectivity radius of $\widehat{\mathfrak S}_l$]\label{prop:injrad}
		For $l>0$, the injectivity radius $\mathrm{InjRad}\left(\widehat{\mathfrak S}_l\right)$ of $\widehat{\mathfrak S}_l$ is equal to $l$.
	\end{proposition}
	\begin{proof}
		Consider a simple closed geodesic curve in $\widehat{\mathfrak S}_l$. If it is contained in a pair of pants, then it exactly matches one of the boundaries and so has length $2l$. Otherwise, we claim that we can find two different pairs of pants in which the curve enters before immediately leaving it through the same boundary component. Indeed, let us consider the projection $c$ of this curve on $\mathfrak T_3$. It is a closed path in $\mathfrak T_3$. Consider $u$ the deepest vertex on $c$ (\textit{i.e.} the vertex with the largest height), and $v$ the furthest point from $u$ in $c$. The vertex that $c$ crosses just before and just after $u$ must be its parent (otherwise $u$ would not be the deepest vertex), and the vertex that $c$ crosses just before and just after $v$ is the unique neighbour of $v$ on the geodesic path from $u$ to $v$ in the tree. It tells us precisely that the two pairs of pants we are looking for are $u$ and $v$.
		Let us take a look at the part $\gamma$ of the curve inside $u$. Without loss of generality, we assume that it enters and leaves $u$ at boundary $\partial u^-$ and refer to Figure \ref{fig:halfpant} and Definition \ref{defi:curves} for the numbering of the shortest paths between boundary components of $P$. Note that $\gamma$ cannot cross only $\gamma_1$, only $\gamma_2$ or none of the $\gamma_i$: otherwise $\gamma$ is contained in the boundary of $P$ and the initial path from which it is extracted can be reduced by not entering $P$, which is not possible since it is a geodesic. For $1 \leq i < j\leq 3$, the distance between $\gamma_i$ and $\gamma_j$ in $u$ is exactly $l$, the length of the half-boundary component between them. Consequently, if $\gamma$ crosses two boundary components in $P$, its length should be at least $l$. If not, then $\gamma$ only crosses $\gamma_3$ and can be decomposed into two paths: one from $\partial P^-$ to $\gamma_3$ and one from $\gamma_3$ to $\partial P^-$. But the shorter curve between $\gamma_3$ and $\partial P^-$ is the unique path which is orthogonal to both of them, (one half of the curve $\widetilde{\gamma}$ shown in red in Figure \ref{fig:halfpant}).	
		By relations for right-angled pentagons (see \cite[Theorem 2.3.4]{Buser}), \begin{align*}l\left({\gamma}\right) &\geq 2 \cosh^{-1}\left(\sinh \frac{l}{2} \ \sinh\cosh^{-1}\left(\frac{\cosh {l} + \cosh^2 {l}}{\sinh^2 {l}}\right)\right)\\
			&= 2 \cosh^{-1}\left(\sqrt{\frac{\cosh^2 {l} + 2\cosh^3 {l} + 1}{\sinh^2 l}}\right)\\ &\geq 2 \cosh^{-1}\left(\sqrt{1 + \cosh{{l}}}\right) = 2\cosh^{-1}\left(\sqrt{2}\cosh\frac{l}{2}\right) \geq {l}.\end{align*}
		The same happens in $v$, and so the curve we have chosen has a length of at least $2l$. So the injectivity radius, which is half the length of the shortest simple closed geodesic, is $l$.
	\end{proof}
	
	\begin{proposition}\label{prop:growthgood}
		There exists a numerical constant $C>0$ such that for $l$ large enough and $R\geq 2\Delta_+$, \[\#\mathcal B_{R+l}\geq C{e^{l}}\# \mathrm{Good}_R,\]
	\end{proposition}
	
\begin{figure}[ht]
		\centering
		\subfloat[0.5\textwidth][$B_h(x_P, l)$ in the first case: in red, the part of the ball in $P$. $\partial P^-$ in green delimits half of the ball, contained in $\overline{\mathcal{D}}(P)$.]
		{
			\centering
			\includegraphics[width=0.35\textwidth]{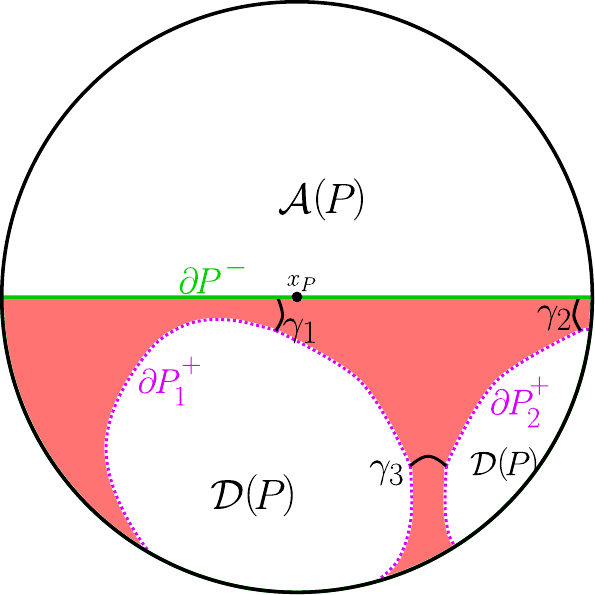}
			\label{fig:ep1}
		}
		\hfill
		\subfloat[0.5\textwidth][$B_h(x_P, l)$ in the second case: in red, the part of the ball in $P$. $\partial P^-$ in green and $\widetilde{\gamma}$ in blue delimitate a quarter of the ball contained in $P$ and the unexplored connected component $b$ of $\mathcal D(P)$.]
		{
			\centering
			\includegraphics[width=0.35\textwidth]{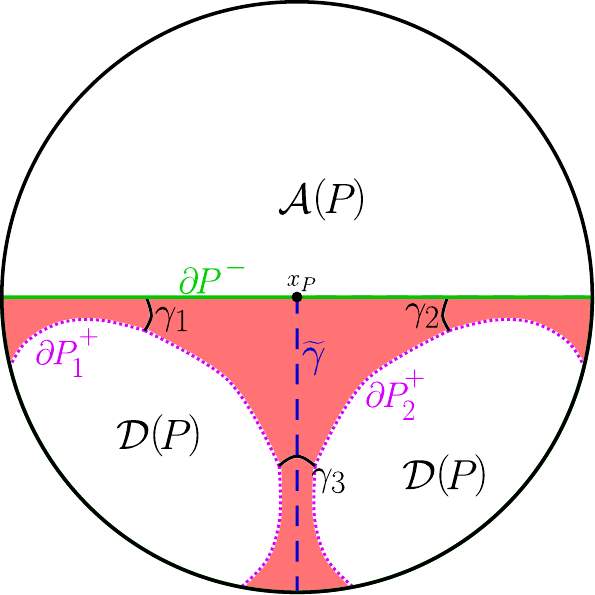}
			\label{fig:ep2}
		}
		\caption{The two possible cases in the proof of Proposition \ref{prop:growthgood}.}
		\label{fig:ep}
	\end{figure}

	\begin{proof}
		Because a pair of pants is of area $2\pi$, it is enough to prove \[\Vol(\mathcal B_{R+l}) \geq C{e^{l}}\# \mathrm{Good}_R,\] up to changing the value of $C$.
		
		The main idea is to find a family $(E_P)_{P \in \Good_R}$ of disjoint subsets of $\mathcal B_{R+l}$, each of which is of volume at least $Ce^{l}$. The fact that $P$ is in $\Good_R$ ensures that there is a point $x_P \in P \cap B_h(\partial \rho^-, R)$ not too far from an unexplored branch in the descendancy of $P$, where the growth rate is large enough. Therefore, $E_P$ is chosen inside this unexplored branch, which is essential to get the disjointness of all the $(E_P)_{P\in\Good_R}$.
		
		Let us describe the sets $E_P$.
		\begin{itemize}
			\item Take $P\in \mathcal S_R, \mathcal D(P)\cap\mathcal B_R = \emptyset$. Fix $x_P \in \partial P^-$ such that $d_h(\partial \rho^-, x_P) \leq R$. The ball $B_h(x_P, l)$ is included in $\mathcal B_{R+l}$. Since $R \geq 2\Delta_+ > 2l$, we are far enough from the boundary of ${\mathfrak S}_l$, and so the injectivity radius at $x_P$ is at least $l$, the one of $\widetilde{\mathfrak S}_l$ (Proposition \ref{prop:injrad}). Therefore, the ball $B_h(x_P, l)$ is isometric to a ball in the hyperbolic plane, so has volume at least $2Ce^l$ if $l$ is large enough and $C > 0$ is small enough. An illustration is given in Figure \ref{fig:ep1}. $\partial P^-$ is one of its diameters, and so half of the volume is contained in $E_P = \overline{\mathcal D}(P) \cap B_h(x_P, l)$.
			\item Take now $P \in \mathrm{Good}_R$ such that $E_P$ has not yet been defined. We fix $x_P = \partial P^-(l/2)$. The ball $B_h(x_P, l-1)$ is included in $\mathcal B_{R+l}$ since $P \in \mathrm{Good}_R$, and $\partial P^-$ is one of its diameter.  The geodesic $\widetilde{\gamma}$ orthogonal to $\partial P^-$ at $x_P$ -- already introduced in the proof of Proposition \ref{prop:injrad}, see Figure \ref{fig:halfpant} -- is also a radius of $B_h(x_P, l-1)$ (see Figure \ref{fig:ep2}). Since $P \in \mathcal S_R$, at least one of the two connected components of $\mathcal D(P)$ does not intersect $\mathcal B_R$. Let $b$ be this connected component. As a consequence, $\partial P^- \cup \widetilde{\gamma}$ delimits a quarter $E_P$ of the ball which is contained in $P \cup b$. As before, the volume of $E_P$ is at least $Ce^l$ if $C$ is small enough and $l$ is large enough.
		\end{itemize}
		
		With such a construction, remark that $P$ is always the youngest ancestor in $\mathcal S_R$ of any pair of pants in $E_P$, and therefore all the $E_P$ are disjoint from each other.
	\end{proof}
	
	Consequently, only the pairs of pants in $\mathcal S_R \setminus \mathrm{Good}_R$ could reduce the growth rate. Rather than finding a way to deal with them, we show that a sufficiently high proportion of pairs of pants in $\mathcal B_R$ are also in $\mathrm{Good}_R$, at least in expectation. This is where the uniformity will be used.
	\begin{proof}[Proof of Theorem \ref{th:balpha}]
		First note that by Proposition \ref{prop:sphere/ball}, $\#\{P\in \mathcal S_R, \mathcal D(P)\cap\mathcal B_R = \emptyset\} \geq \#(\mathcal B_R\setminus \mathcal S_R)$, and so \begin{align*}2\#\Good_R &\geq \#\left \lbrace P\in \mathcal B_R\setminus \mathcal S_R \text{ or } P \in \mathcal S_R,\  \partial P^-\left(\left[\frac{l}{2}-1; \frac{l}{2}+ 1\right]\right) \cap B_h(\partial \rho^-, R)\neq \emptyset\right\rbrace \\&\geq \#\left \lbrace P \in \mathcal B_R,\  \partial P^-\left(\left[\frac{l}{2}-1; \frac{l}{2}+ 1\right]\right) \cap B_h(\partial \rho^-, R)\neq \emptyset\right\rbrace,
		\end{align*}
		where the factor $2$ accounts for pairs of pants which respect both conditions in $\Good_R$ and are likely to be counted twice in the right-hand side.
		For a point $P$, we write $\Pi_P$ for the point on $\partial P^-$ closest to $\partial \rho^-$ (which is the first point to be "seen" by the spheres $\left(\mathcal S_R\right)_R$), and $t_{P} \in [0, 2l[$ the parameter such that $\Pi_P= \partial P^- (t_{P})$. Then, the previous inequality implies that
		\begin{align*}\E\#\mathrm{Good}_R &\geq \frac{1}{2}\E\#\left \lbrace P \in \mathcal B_R,\  \partial P^-\left(\left[\frac{l}{2}-1; \frac{l}{2}+ 1\right]\right) \cap B_h(\partial \rho^-, R)\neq \emptyset\right\rbrace \\&\geq \frac{1}{2} \bigop{\sum}{P \in V^\mathfrak B}{} \E\left[\boldsymbol{1}_{P\in\mathcal B_R}\boldsymbol{1}_{t_{P} \in \left[\frac{l}{2}-1; \frac{l}{2}+ 1\right]}\right].
		\end{align*}
		Let us fix $P \in V^\mathfrak B$. On the one hand, the event $P \in \mathcal B_R$ (\textit{i.e.} the event $d_h(\partial \rho^-, P) \leq R$) is measurable with respect to the $\sigma$-algebra $\mathcal F$ generated by the twists and lengths inside $\mathcal A(P)$ and all the lengths of boundary components of $\mathcal A(P)$. In particular, it does not depend on the twist between $P$ and its parent. On the other hand, because of the choice of parametrization, $t_P$ is a uniform random variable independent of $\mathcal F$: it is the sum of an $\mathcal F$-measurable random variable -- which is the position of $\Pi_P$ on the boundary of the parent of $P$ -- and a uniform random twist between $P$ and its parent which is independent of $\mathcal F$. As a consequence, \begin{align*}\bigop{\sum}{P \in V^\mathfrak B}{} \E\left[\boldsymbol{1}_{P\in\mathcal B_R}\boldsymbol{1}_{\Pi_P \in \partial P^-\left(\left[\frac{l}{2}-1; \frac{l}{2}+1\right]\right)}\right] &= \bigop{\sum}{P \in V^\mathfrak B}{} \E\boldsymbol{1}_{P\in\mathcal B_R}\E\boldsymbol{1}_{t_P \in \left[\frac{l}{2}-1; \frac{l}{2}+ 1\right]}\\&= \frac{2}{2l} \bigop{\sum}{P \in V^\mathfrak B}{} \E\boldsymbol{1}_{P\in\mathcal B_R},\end{align*} so \[\E\#\mathrm{Good}_R \geq \frac{1}{2l} \E N_R.\]
		
		From Proposition \ref{prop:growthgood}, and using Proposition \ref{prop:sphere/ball}, we therefore obtain that for $l$ large enough, \[\forall R\geq 2\Delta_+,\quad \E N_{R+l}\geq \frac{C}{4l} {e^{l}}\E N_R.\]
		Iterating this $\lfloor \ln l \rfloor$ times from $R = 2\Delta_+$, we get that \[\frac{\ln\E N_{2\Delta_+ +l\lfloor \ln l \rfloor}}{2\Delta_+ +l\lfloor \ln l \rfloor} \geq \frac{l\lfloor \ln l  \rfloor - o(l \lfloor\ln l\rfloor)}{2\Delta_+ +l\lfloor \ln l \rfloor}.\]
		
		Since $\Delta_+ = o(l \lfloor\ln l\rfloor)$ according to Proposition \ref{prop:Ddelta}, and thanks to Corollary \ref{coro:convexp}, we deduce that for $\varepsilon >0$ and $l$ large enough, \[\Sigma_{l \ln l} \geq 1-\varepsilon, \qquad \left \lvert \alpha - \Sigma_{l\ln l}\right\rvert \leq \varepsilon.\] As a consequence, \[\liminf_{l \to \infty} \alpha \geq 1,\] as expected. 
	\end{proof}
	
	Note that the only moment we specifically use the fact that twists are uniform is when we show that $\E\mathrm{\#Good}_R$ is "large enough". The author believes this fact should remain true for a wide range of twist law: indeed, the sequence of $(t_P)_{P \in V^\mathfrak B}$ should be approximated by the steps of a branching random walk over $\R/l\Z$, and consequently, the number of good pairs of pants should approximately be controlled by a kind of Cesaro mean of this random walk, which roughly behave uniformly. Nevertheless, questions of independence and conditioning make things harder to prove in the general setting.
	
	\section{Back to \texorpdfstring{$\mathfrak S_{\nu, g}$}{S\_\{ν,g\}}}
	\label{Surface}
	
	For this section, we recall the following notations :
	\begin{itemize}
		\item $u_g = o(v_g)$ means that for all $\varepsilon > 0$, for $g$ large enough $u_g \leq \varepsilon v_g$.  We also use the notation $u_g \ll v_g$.
		\item $u_g = O(v_g)$ means that there exists $C>0$ independent of $g$ such that $u_g \leq C v_g$.
	\end{itemize}
	
	Let us recall that $\mathfrak S_{\nu, g}$ is built by taking iid (length, twists)-weights with law $\nu = \nu_g$ on a 3-regular graph with $2g-2$ vertices built from the configuration model. More precisely, start from $2g-2$ vertices to each of which 3 half-edges are attached, and pair the half-edges uniformly at random to get the edges of the graph. Concretely, it corresponds to pair uniformly at random the boundary components of the $2g-2$ pairs of pants used to build the surface, and glue the paired boundary components together with lengths and twists being independent from one gluing to another, with distribution $\nu$. Note that the pairing can be realized algorithmically: fix arbitrarily any unpaired half-edge, and then choose uniformly at random among all other unpaired half-edges the one to be its match for creating an edge of the graph, and keep doing it until all the half-edges are paired.
	
	This section is dedicated to proving Theorem \ref{th:Principal3}. More precisely, we will prove the following reformulation: up to taking $\eta/2$ instead of $\eta$ -- or $\sqrt{u_g}$ instead of $u_g$ for the uniformly bounded case -- in the following statement, it is enough to obtain Theorem \ref{th:Principal3}. Note that the bounds on $\alpha$ are direct consequence of Lemma \ref{lem:majcroiss}, as remarked at the beginning of Section \ref{behavalpha}.
	\begin{theorem}\label{th:Simp}
		Suppose that \[\supp\left(\nu_l\right) \subset [a_g^-, a_g^+], \quad \text{with } \frac{1}{a_g^-} = o(\ln \ln g), a_g^+ = o (\ln \ln g).\] 
		Take any $\eta > 0$, and set $\zeta_g = \ln ^{3/4+\eta}g$. Then asymptotically almost surely \[\left\lvert {\diam \mathfrak S_{\nu, g} - \frac{1}{\alpha_{\nu}}\ln g}\right \rvert \leq \zeta_g,\] with $\alpha_\nu$ defined in Corollary \ref{coro:convexp}.
		Moreover, when $\supp{\nu}$ is uniformly bounded in $g$, we can take $\zeta_g = u_g\ln^{3/4}g$ for any $(u_g)_{g\geq 0}$ such that $u_g \limitginf{}{} \infty$. 
	\end{theorem}
	
	The proof of Theorem \ref{th:Simp} is based on a precise analysis of an exploration process of the graph and of the associated surface. In the Section \ref{Algo}, we describe the exploration and give technical lemmas to compare $\mathfrak{S}_{\nu, g}$ and $\mathfrak{S}_\nu$ defined in Section \ref{TLS}. In Section \ref{Prooffin}, we use these lemmas to prove Theorem \ref{th:Simp}.
	
	\subsection{Exploration of the graph}\label{Algo}
	
	As in \cite{BCP_mindiam}, the proof is mainly based on a detailed study of a kind of breadth-first exploration around a fixed vertex $\rho$ of the graph used to build $\mathfrak{S}_{\nu, g}$, taking into account the hyperbolic metric. Starting from $\rho$, we pair sequentially the edges to discover its neighbourhood. The exact description of the algorithm to select the half-edge to be paired at each step is postponed to Definition \ref{defi:algo}. We introduce the following vocabulary.
	\begin{definition}
		If the connected component of $\rho$ is entirely discovered before we decide to stop the algorithm, we say that the exploration \textbf{stops prematurely} (in particular, this means that the surface is disconnected). At each pairing, two things can happen: either we add a new vertex to the neighbourhood we have discovered (by pairing the current half-edges to one belonging to a pair of pants which has not been seen before), or we pair two pairs of pants (or one with itself) already discovered. In this last case, we say that a \textbf{bad step} occurred, following the terminology of \cite{BCP_mindiam}.
	\end{definition}
	
	In the proof of Theorem \ref{th:Simp}, we study what happens for balls when they contain approximately $\sqrt{g}$ pairs of pants. By a kind of birthday paradox argument, we get that it is around this value that a transition happens: before containing $\sqrt{g}$ pairs of pants, two balls will not intersect each other, while they become very likely to do so after. Therefore, we study these two "regimes". First, we try to understand what happens from a metric point of view if we stop the exploration as soon as slightly more than $\sqrt{g}$ vertices are found, and justify that the radius of the explored part is not much larger than the radius of a ball in $\mathfrak S_\nu$ containing the same number of pairs of pants. Then we also have to argue that bad steps do not increase significantly the growth rate (by creating long "shortcuts" thanks to cycles in the graph and reducing distances between points), at least if we stop when we find slightly less than $\sqrt{g}$ vertices. Note that in \cite{BCP_mindiam} and \cite{Bollobas}, this last point was not a difficulty, since a comparison with the universal cover automatically gives a sufficient lower bound on the diameter. These facts are gathered in the following lemmas, which are the aim of this subsection. In what follows, "the exploration" refers to the one described by Definition \ref{defi:algo}.
	
	\begin{lemma}[{Bounds on the number of bad steps}]\label{lem:borneserreurs}
		Let $0< \beta < 1/2$, and $k > \frac{2}{1-2\beta}$. Consider the exploration around $\rho$ until $\sqrt{g}\ln g$ vertices are found. Then, with probability $1-o\left(g^{-2}\right)$, either the exploration stops prematurely, or \begin{itemize}
			\item less than $k$ bad steps are made during the first $g^\beta$ steps,
			\item less than $\ln^{3/4} g$ bad steps are made before finding ${g}^{\frac{1}{2}-\frac{1}{\ln^{3/4}g}}$ vertices,
			\item and less than $\ln^3 g$ bad steps are made during the whole exploration
		\end{itemize} at the same time.
	\end{lemma}
	
	\begin{lemma}[{Upper bound for the radius of the exploration}]\label{lem:detmaj}
		Suppose that the exploration around $\rho$ does not stop prematurely, that less than $k$ bad steps have been made during the first $g^\beta$ first steps, and less than $\ln^3 g$ bad steps have been made until we find $\sqrt{g}\ln g$ vertices.
		Then with probability $1-o\left(g^{-2}\right)$, any point discovered at this stage is at most at distance $\frac{1}{2\alpha}\ln g + \zeta_g$ of $\rho$ in $\mathfrak S_{\nu, g}$, where $\zeta_g$ is introduced in Theorem \ref{th:Simp}. 
	\end{lemma}
	
	\begin{lemma}[{Lower bound for the radius of the exploration}]\label{lem:detmin}
		Suppose that the exploration does not stop prematurely before we find $g^{\frac{1}{2}-\frac{1}{\ln^{3/4}g}}$ vertices, and that less than $\ln^{3/4} g$ bad steps are made. Then with probability $1-o\left(g^{-2}\right)$, any unpaired half-edge is at least $\frac{1}{2\alpha}\ln g - \zeta_g$ away from $\rho$, where $\zeta_g$ is introduced in Theorem \ref{th:Simp}. In particular, any point at distance less than $\frac{1}{2\alpha}\ln g - \zeta_g$ from $\rho$ in $\mathfrak S_{\nu, g}$ is contained in the discovered neighbourhood. 
	\end{lemma}
	
	Let us first begin with the proof of Lemma \ref{lem:borneserreurs}, which does not depend on the exact description of the exploration. We only impose that our algorithm respect the construction of the configuration model given in Definition \ref{defi:config}: once you have picked a boundary in the explored part, the boundary you pair it with is uniform among all available boundaries.
	\begin{proof}[Proof of Lemma \ref{lem:borneserreurs}]
		It is a similar lemma as the one used in the proof of \cite[Proposition 3]{BCP_mindiam}, itself inspired by what has been done in \cite{Bollobas} for graphs. In any case, at most $3\sqrt{g}\ln g$ steps are made. At step $i$, the number of discovered edges still unpaired is bounded from above by $3+2i$, and the number of undiscovered edges is bounded by $3g-2i$. As a consequence, for $1\leq i \leq 3\sqrt{g}{\ln g}$, the probability of step $i$ to be a bad step is bounded from above by $C\frac{i}{g}$ for some $C>0$ which does not depend on $g$ and therefore, \[\p\left(\text{more than $k$ bad steps before time $g^\beta$}\right) \leq \binom{g^\beta}{k}C^kg^{k(\beta-1)} \leq \frac{C^k}{k!}g^{k(2\beta-1)}.\] After $t$ steps, $2t$ half-edges have been paired, and so at least $2t/3$ vertices have been found, because each vertex has three half-edges. So for any $\widetilde{\beta}$, we need at most $\frac{3}{2}g^{\widetilde{\beta}} < 2g^{\widetilde{\beta}}$ steps to find $g^{\widetilde{\beta}}$ vertices. As before, using Stirling equivalent, for any $\widetilde{\beta}$ such that $g^{\widetilde{\beta}} \leq \sqrt{g}\ln g$, for any $b^2 \ll g^{\widetilde{\beta}}$,
		\begin{align*}\p\left(\text{more than $b$ bad steps before finding $g^{\widetilde{\beta}}$ vertices}\right) &\leq \binom{2g^{\widetilde{\beta}}}{b}(2C)^bg^{b(\widetilde{\beta}-1)}\\&\leq \frac{1}{\sqrt{b}} \left(\frac{4Ceg^{2\widetilde{\beta}-1}}{b}\right)^b.\end{align*}
		Take $\widetilde \beta = \frac{1}{2}-\frac{1}{\ln^{3/4}g}$, $b = \ln^{3/4} g$ or $\widetilde \beta = \frac{1}{2} + \frac{\ln\ln g}{\ln g}$, $b = \ln^{3} g$ : these upper bounds are both  $o\left(g^{-2}\right)$.
	\end{proof}
	
	In \cite{BCP_mindiam} they choose $\beta = \frac{1}{2}-\varepsilon$, but it seems to us that any $0<\beta<1/2$ will suffice for the proof. The choice of $\sqrt{g}\ln g$ is more crucial: a value larger than $\sqrt{g}$ is needed to ensure that two such explorations around two different vertices merge with a sufficiently high probability, but cannot be too large so as not to degrade the growth rate obtained. On the contrary, the value just below $\sqrt{g}$ for the second exploration will guarantee that two explorations around two different vertices will not merge. It will be essential for the proof of Theorem \ref{th:Simp}.
	
	This lemma expresses the fact that when we discover the graph by a pseudo breadth-first search, we have relatively few bad steps at the beginning. This leads to the existence of many subtrees, and thus to a growth that is roughly the same as that of $\mathfrak S_\nu$, as demonstrated by Lemmas \ref{lem:detmaj} and \ref{lem:detmin}.
	
	Let us now give a precise description of the algorithm of exploration we use in this section.
	
	\begin{definition}\label{defi:algo}
		We use the following \textbf{algorithm for the exploration}. 
		Denote by $(G_i, w_i)$ the neighbourhood of $\rho$ discovered after $i$ steps of the algorithm, where $w_i$ corresponds to the weights put on the edges of $G_i$. Let $E_i$ be the set of unpaired half-edges of $G_i$ (\textit{i.e.} the set of boundary components of the corresponding surface). For example, $G_0$ is the graph with only one vertex $\rho$, $E_0$ is the set of the three half-edges of $\rho$ and $w_0 = \emptyset$ (there is no edge in $G_0$, only half-edges). At the next step,
		\begin{itemize}
			\item if $E_i = \emptyset$, we say that the exploration stops prematurely and the algorithm stops.
			\item Otherwise, consider the topological surface described by the graph $G_i$. To get a hyperbolic surface (\textit{i.e.} to define a hyperbolic metric on it), we should use the weights $w_i$, but also prescribe length for all the boundary components $e \in E_i$. Let $\mathfrak M(G_i, w_i)$ be the set of all hyperbolic surfaces that can be obtained for fixed $G_i$ and $w_i$. For any unpaired boundary component $e \in E_i$, let \[\bm{{d_{+}^{G_i}(e)}} = \max_{S\in \mathfrak M(G_i, w_i)} d_h(\rho \subset S, e \subset S).\] Note that $d_{+}^{G_i}(e)$ implicitly depends on $w_i$, but we do not make this dependence explicit in the notations for simplicity. Select the next half-edge to be paired $e^*_i$ among the minimizers of $d_{+}^{G_i}$. Finally, choose the weight of the new edge according to $\nu$. 
		\end{itemize}
		More generally, for any graph $(G, w)$ with weights on its edges and unpaired half-edges, and $x$ in the associated topological surface, we write \[\bm{{d_{+}^{G}(x)}} = \max_{S\in \mathfrak M(G, w)} d_h(\rho \subset S, x \subset S).\]
	\end{definition}
	
	Concretely, at each step, we consider the best upper bound we can have on the distances from $\rho$ to the boundary components, taking into account the weights we have already discovered, and select the edges for which the upper bound is the smallest. Note that this would not always be the closest edge for the graph distance, and that for any $i\leq j$ and $x \in G_i$, \[d^{\mathfrak S_{\nu, g}}(\rho, x) \leq d_+^{G_j}(x) \leq d_+^{G_i}(x),\]
	because by definition, $\mathfrak S_{\nu, g} \subset \mathfrak M(G_j, w_j) \subset \mathfrak M(G_i, w_i)$.
	
	Before going deeper into the details, let us state the following technical lemma, which explains why our algorithm can indeed be considered as a kind of breadth-first exploration. It will be particularly useful for both Lemma $\ref{lem:detmaj}$ and \ref{lem:detmin}.
	
	\begin{lemma}\label{lem:breadth}
		Suppose that we stop the exploration at step $i$ and find a point $x$ such that $d_{+}^{G_i}(x) = R$. Then any unpaired boundary component $e$ in $G_i$ is such that $d_h^{\mathfrak S_{\nu, g}}(\rho, e) \geq R-2\Delta_+$. In particular, if $d_+^{G_i}(e) \leq R$, then any edge paired before $e$ is such that $d_+^{G_i}(e) \leq R+ 2\Delta_+$.
	\end{lemma}
	
\begin{figure}[ht]
		\centering
		\includegraphics[width=0.35\textwidth]{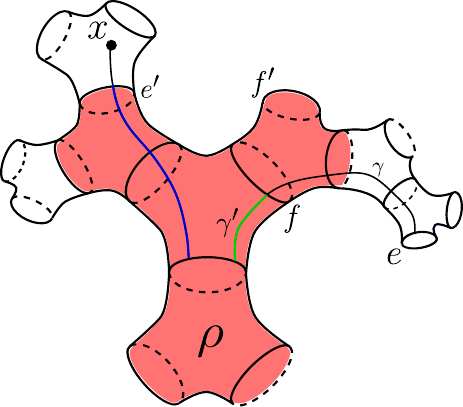}
		\caption{$G_i$ and the path $\gamma$. Note that $\gamma$ might not be entirely drawn in $G_i$ (if it uses gluings after the step $i$). In red, the graph $G_{e'}$, and in green, $\widetilde{\gamma}$ whose length remains unchanged in the next steps because it does not cross pairs of pants with unglued boundary components. Because $d_{+}^{G_i}(x) = R$, we have that $d_{+}^{G_{e'}}(e') \geq d_{+}^{G_i}(e') \geq R-\Delta_+$, and so $\ell(\widetilde{\gamma}) \geq d_{+}^{G_{e'}}(f) \geq d_{+}^{G_{e'}}(f')-\Delta_+ \geq R-2\Delta_+$.}
		\label{fig:lem37}
	\end{figure}

	\begin{proof}
		Consider in $\mathfrak S_{\nu, g}$ a path $\gamma$ from $\rho$ to a boundary component $e$ which was unglued in $G_i$. The aim is to prove that $\ell(\gamma) \geq R-\Delta_+$. We use Figure \ref{fig:lem37} to illustrate the proof.
		
		Let $P$ be the pair of pants in $G_i$ containing the point $x$ such that $d_{+}^{G_i}(x) = R$. Let $e'$ be a glued boundary component of $P$ such that $d_{+}^{G_i}(e') \leq R$ (any path from $\rho$ to $x$ of length $R$ has to cross such a boundary of $P$). Then, $d_{+}^{G_i}(e') \geq R - \Delta_+$ because $d_{+}^{G_i}(x) = R$. Let $G_{e'}$ be the part of $G_i$ we had already explored just before gluing $e'$. It is in particular possible that $P$ does not belong yet to $G_{e'}$, but $e' \in G_{e'}$ as part of the pair of pants which is going to be paired with $P$. In that graph, $d_{+}^{G_{e'}}(e') \geq R-\Delta_+$ as per the remark just after Definition \ref{defi:algo}. But $e'$ is the edge selected to be paired, so any other unpaired half-edge $f$ verifies $d_{+}^{G_{e'}}(f) \geq R-\Delta_+$. Consider $\gamma'$ the part of $\gamma$ starting from $\rho$ and stopping just before it enters a pair of pants $\widetilde{P}$ with an unglued boundary component. In particular, $\gamma'$ might be empty (if $\rho$ has an unglued boundary) or reduced to a point (if the first pair of pants crossed by gamma has an unglued boundary in $G_{e'}$). Let $f$ be the boundary components where $\gamma'$ ends. $\widetilde{P}$ contains an unglued boundary component $f'$, which is such that $d_{+}^{G_{e'}}(f') \geq R-\Delta_+$, and so $d_{+}^{G_{e'}}(f) \geq R-2\Delta_+$. In particular, $l_{G_e'}(\gamma') \geq R-2\Delta_+$. But $\gamma'$ does not cross any pair of pants with unpaired boundary components, so its length will not be changed during the end of the exploration. In particular, $\ell(\gamma)\geq R-2\Delta_+$.
		
		For the second affirmation, suppose by contraposition that we can find an edge $e'$ paired before $e$ such that $d_+^{G_i}(e') > R+ 2\Delta_+$. Let $j \leq i$: then $d_+^{G_{j}}(e') > R+ \Delta_+$, and apply the first statement to $G_{j}, x \in {e'}$ to get that any boundary component $f$ paired after $e'$ but before $e$ should verify $d_+^{G_j}(f) > R$. This applies in particular to $f=e$ itself.
	\end{proof}
	
	Note that a classical breadth-first search exploration would not be possible in our case, and that is why we choose this modified version. Indeed, we cannot assign a length to a boundary component before knowing which one is paired with which one, and so we cannot know the exact geometry of the discovered neighbourhood.
	
	In the following, we give separately the tools needed to prove the upper bound and those for the lower bound. Both will be based on coupling lemmas.

	Before we go into the details, it is worth noting where the bounds on the support of $\nu_l$ appear. In what follows -- more especially in Lemmas \ref{lem:detmaj} and \ref{lem:detmin} -- we often need to compare some quantities expressed in terms of $\Delta_+, \delta_-, \alpha$ with $\zeta_g$ introduced in Theorem \ref{th:Simp}. Since we already have estimates of $\Delta_+, \delta_-, \alpha$ in terms of the support of $\nu_l$ -- Property \ref{prop:Ddelta} and Lemma \ref{lem:majcroiss} -- the hypothesis on the support of $\nu_l$ is here to let us make these comparisons.
	In particular, when $\mathrm{Supp}(\nu_l)$ is uniformly bounded, we see that it is also the case of $\Delta_+, \delta_-, \alpha$, and so taking $\zeta_g = u_g\ln^{3/4}g$ for a sequence $(u_g)_{g\geq 0}$ such that $u_g \limitginf{}{} \infty$ will be sufficient.
	
	\subsubsection{Tools for the upper bound}
	The aim of this section is to obtain Lemma \ref{lem:detmaj}, which gives an upper bound on the distances from $\rho$ of the $\sqrt{g}\ln g$ first vertices discovered during the exploration around $\rho$.
	
	We recall that the notation $^{\mathtt C}{}{T}_i$ introduced in Definition \ref{defi:childtree}, and appears here to preserve all the information contained in $\mathcal F_{T_i}$, to capture the whole geometry of $T_i$.
	We will use the following coupling:	
	\begin{lemma}\label{lem:couplage}
		Take a weighted graph $(G, w^G)$, all of whose weights are chosen independently with distribution $\nu$. Consider a family $I$ of rooted subgraphs $(T_i, \rho_i, w^{^{\mathtt C}{}{T}_i})_{i \in I}$ of $G$ pairwise at graph distance 2, all isomorphic to subtrees of $(\mathfrak B, \rho)$. Suppose in addition that any two leaves of $T_i$ are also at graph distance 2 -- which amounts to saying that $^{\mathtt C}{}{T}_i$ is also a tree in $G$. This family can be coupled with a family $(S_i, \rho_i, w^{S_i})$ of independent copies of $(\mathfrak S_\nu, \rho, w^{\mathfrak S_\nu})$, such that $(T_i, \rho_i, w^{^{\mathtt C}{}{T}_i})$ is isometric to a subtree of $(S_i, \rho_i, w^{S_i})$.
	\end{lemma}
	\begin{proof}
		Consider the weighted graph $(T_i, \rho_i, w^{^{\mathtt C}{}{T}_i})$. It is a finite binary tree, which can be completed in an infinite binary tree, all of whose unknown weights are chosen independently, with law $\nu$. This gives a weighted tree/surface $S_i$, distributed like $\mathfrak S_\nu$. Because the $T_i$ are at graph distance 2, all the weights are independent of one another, and so are the $S_i$.
	\end{proof}
	
	Note that in particular, for every $R_i \leq d_h(\partial \rho_i^-, S_i\setminus T_i)$, \textit{i.e.} $\mathcal B_{R_i}^{S_i}(\rho_i) \subset T_i$, the numbers of vertices/pairs of pants at hyperbolic distance less than $R_i$ from $\partial \rho_i^-$ in $T_i$ are equal to $\left(N_{R_i}(\rho_i)\right)_i$, and are distributed like independent copies of $N_{R_i}$.
	
	Thanks to the coupling argument, we can now prove Lemma \ref{lem:detmaj}, analogous in our model to \cite[Lemma 4]{BCP_mindiam}, which expresses the fact that $\mathfrak S_{\nu, g}$ does not have a significantly worse behaviour than $\mathfrak S_\nu$.
	
	\begin{proof}[Proof of Lemma \ref{lem:detmaj}]
		Let us start with a first phase of the exploration -- $\lfloor g^\beta \rfloor$ steps. We write $G$ for the graph obtained at the end of this phase, and $R$ the maximal hyperbolic distance from $\rho$ possibly reached in $G$: \[ R = \max_{x \in G} d_+^{G}(x).\] Since the exploration is not stopped prematurely and there are at most $k$ bad steps, we can find a vertex $v_0$ which is at graph distance $k$ from the root $\rho$, whose descendancy in $(G, \rho)$ is a tree (no bad steps). Let $h = h_g = \left \lceil \log_2 \ln g \right \rceil$ and consider the family $(\rho_i)_{1 \leq i \leq 2^h}$ of the $2^h$ descendants at the $h^{\mathrm{th}}$ generation of $v_0$. They all have a tree-like descendancy in $(G, \rho)$. Let $(T_i)_{1 \leq i \leq 2^h}$ be these tree-like descendancies where we erased the nodes with unpaired half-edges (in particular, any weight on half-edges adjacent to $T_i$ has been determined while building $G$). The $(T_i)_i$ are at a graph distance 2 from one another because a path from one to another has to go through an erased node or through $v_0$, and so they respect hypotheses of Lemma \ref{lem:couplage}: we then use a coupling $T_i \hookrightarrow S_i$ as in Lemma \ref{lem:couplage} to estimate $R$.
		
		According to Lemma \ref{lem:breadth}, if we find an unpaired boundary component $e$ in $g$ such that $d_+^G(e) \leq R' - 2\Delta_+$ for some $R'$, then $R\leq R'$. Since $G$ contains less than $g^\beta$ vertices if we find $(i, \widetilde{R})$ such that $N_{\widetilde R}(\rho_i) \geq g^\beta$, then $\widetilde R \geq d_h(\partial \rho_i^-, S_i\setminus T_i)$, as remarked after Lemma \ref{lem:couplage}. This gives an unglued boundary component at distance at most  \begin{align*} d_+^{G}(\rho_i) + \Delta_+ +  d_h(\partial \rho_i^-, S_i\setminus T_i) + \Delta_+
			&\leq R' = (k+h+3) \Delta_+ + \widetilde{R}
		\end{align*}
		from $\rho$, where the last $\Delta_+$ in the first line accounts for the fact that there is a node between $T_i$ and the boundary of $G$. Consequently \[R \leq (k+h+5) \Delta_+ + \widetilde{R}.\]
		Take $\widetilde{R} = \frac{\beta}{\alpha} \ln g + \zeta_g$, and $\varepsilon = \alpha - \frac{\beta \ln g}{\widetilde{R}} = \frac{\alpha^{2} \zeta_g}{\beta\ln g} + o\left(\frac{\alpha^{2} \zeta_g}{\ln g}\right)$  so that $0<\varepsilon<1$, $\varepsilon \limit{}{g \to \infty} 0$ because \textbf{$\alpha \zeta_g \ll \ln g$}.
		Since $\zeta_g \gg \frac{e^{c\Delta_+}}{\alpha^{2}}\ln^{3/4}g$, then \textbf{$\widetilde{R} \gg \frac{e^{c\Delta_+}}{\varepsilon^4}$}, so that Theorem \ref{th:concentration} applies for $g$ large enough. This gives that 
		\begin{align*}\p\left (N_{\widetilde{R}} \leq g^\beta \right)
			&\leq 2 \exp\left(-e^{-C\Delta_+} \varepsilon^2\widetilde{R}\right)
			\\& \leq 2 \exp\left(-\frac{\alpha^3\zeta_g^2}{e^{C\Delta_+}\beta \ln g} + o\left(\frac{\alpha^3\zeta_g^2}{e^{C\Delta_+}\ln g}\right)\right)\\ &= \exp\left(-\frac{\alpha^3\zeta_g^2}{e^{C\Delta_+}\beta \ln g} + o\left(\frac{\alpha^3\zeta_g^2}{e^{C\Delta_+}\ln g}\right)\right),\end{align*}
		because \textbf{${\alpha^3\zeta_g^2}\gg {e^{C\Delta_+}\ln g}$}. In this way, by independence of the $T_i$, \begin{align*}\p\left(\forall 1\leq i \leq 2^h, N_{\frac{\beta}{\alpha} \ln g + \zeta_g}(\rho_i) \leq g^\beta\right) &\leq \exp\left(2^h \left(-\frac{\alpha^3\zeta_g^2}{e^{C\Delta_+}\beta \ln g} + o\left(\frac{\alpha^3\zeta_g^2}{e^{C\Delta_+}\ln g}\right)\right)\right)\\&\leq g^{-\frac{\alpha^3\zeta_g^2}{e^{C\Delta_+}\beta \ln g} + o\left(\frac{\alpha^3\zeta_g^2}{e^{C\Delta_+}\ln g}\right)} = o\left(g^{-2}\right),\end{align*} because \textbf{${\alpha^3\zeta_g^2}\gg {e^{C\Delta_+}\ln g}$}.
		
		Thus, with probability $1-o\left(g^{-2}\right)$, one of the $N_{\frac{\beta}{\alpha} \ln g + \zeta_g}(\rho_i)$ verifies $N_{\frac{\beta}{\alpha} \ln g + \zeta_g}(\rho_i) > g^\beta$, and so, \[R \leq (k+5+h)\Delta_++ \frac{\beta}{\alpha} \ln g + \zeta_g\leq\frac{\beta}{\alpha} \ln g + 2\zeta_g,\] because \textbf{$h\Delta_+\ll \zeta_g$}.
		
		Now, take a look at the second phase of exploration, and note $G'$ the resulting graph, $R'$ the maximal hyperbolic distance from $\rho$ in the surface $G'$: \[ R' = \max_{x \in G'} d_+^{G'}(x).\]
		
		At the end of the first step, there are at least $g^{\beta}-3k$ half-edges discovered but still not paired. Indeed, at each step, the number of half-edges discovered but not paired yet increases by 1 (gains 2, loses 1) unless it was a bad step (gains 0, loses 2). Let us continue the exploration until all of these half-edges are paired (let us say that we explore up to radius $R+2\Delta_+$). Because less than $\ln^3 g$ errors are made, we observe a family $I_2$ of at least $\left \lfloor g^{\beta} - 3\ln^3 g\right \rfloor $ pairs of pants at graph distance $2$, whose descendancies must be disjoint trees in $(G_2, \rho)$. As before, we can build a coupling $(T_i \hookrightarrow S_i)_{i \in I_2}$ as per Lemma \ref{lem:couplage}.
		We claim that, as before, with probability $1-o\left(g^{-2}\right)$, among the $\left \lfloor g^{\beta} - 3\ln^3 g \right \rfloor$ independent copies of $\mathfrak S_\nu$, $\left \lceil \frac{g^{\beta}}{\ln g}\right \rceil$ are such that \[N_{\frac{1-2\beta}{2\alpha}\ln g + \zeta_g}(\rho_i) > g^{1/2-\beta}\ln^2 g.\] If so, as before, because $G'$ contains $\sqrt{g}\ln g$ vertices, at least one of the $T_i$ is such that $d_h(\partial \rho_i^-, S_i\setminus T_i) < \frac{1-2\beta}{2\alpha}\ln g + \zeta_g$  with probability $1-o\left(g^{-2}\right)$. Consequently, with probability $1-o\left(g^{-2}\right)$
		\begin{align*}
			R' &\leq R + 3\Delta_+ + \frac{1-2\beta}{2\alpha}\ln g + \zeta_g + \Delta_+ 
			\\&\leq \frac{1}{2\alpha}\ln g + 4\zeta_g.
		\end{align*}
		
		It only remains to prove the claim to conclude. 
		Take $R = \frac{(1 - 2\beta)\ln g}{2\alpha}  + \zeta_g$, and $\varepsilon = \alpha - \frac{(1 - 2\beta)\ln g + 4\ln \ln g}{2R}= \frac{2\alpha^{2} \zeta_g}{(1-2\beta)\ln g} + o\left(\frac{\alpha^{2} \zeta_g}{\ln g}\right)$ so that $0<\varepsilon<1$, $\varepsilon \limitginf{}{} 0$ because {$\ln \ln g \ll \alpha\zeta_g\ll \ln g$}.
		As before, Theorem \ref{th:concentration} applies and so	
		\begin{align*}\p\left (N_R \leq g^{(\frac{1}{2}-\beta)}\ln^2g \right)
			&\leq \exp\left(-\frac{2\alpha^3\zeta_g^2}{e^{C\Delta_+}(1-2\beta)\ln g} + o\left(\frac{\alpha^3\zeta_g^2}{e^{C\Delta_+}\ln g}\right)\right).\end{align*}
		
		As a consequence, the probability that among the $\lfloor g^\beta- 3\ln^3g \rfloor$ trees, at least $\lfloor g^\beta- 3\ln^3g \rfloor - \lceil \frac{g^{\beta}}{\ln g} \rceil$ verify $N_R(\rho_i) < g^{(\frac{1}{2}-\beta)}\ln^2g$ is bounded above by:
		\begin{multline*}\left(\star\right) = \binom{\left\lfloor g^\beta-3\ln^3 g \right\rfloor }{\left\lfloor g^\beta- 3\ln^3g \right\rfloor - \left\lceil \frac{g^{\beta}}{\ln g} \right\rceil}\\\times \exp\left(\left(-\frac{2\alpha^3\zeta_g^2}{e^{C\Delta_+}(1-2\beta)\ln g} + o\left(\frac{\alpha^3\zeta_g^2}{e^{C\Delta_+}\ln g}\right)\right)\left(\left\lfloor g^\beta- 3\ln^3g \right\rfloor - \left\lceil \frac{g^{\beta}}{\ln g} \right\rceil\right)\right).\end{multline*}
		
		Let us estimate the binomial coefficient, using Stirling formula: when $k = o(n)$, then \[\binom{n}{k} = \left(\frac{n e}{k} \right)^k \cdot (2\pi k)^{-1/2} \cdot \exp\left(- \frac{k^2}{2n}(1 + o(1))\right)\] and so 
		\begin{align*}
			\binom{\left\lfloor g^\beta-3\ln^3 g \right\rfloor }{\left\lfloor g^\beta- 3\ln^3g \right\rfloor - \left\lceil \frac{g^{\beta}}{\ln g} \right\rceil} &= \binom{\left\lfloor g^\beta-3\ln^3 g \right\rfloor }{\left\lceil \frac{g^{\beta}}{\ln g} \right\rceil} \\&\leq \binom{\left\lfloor g^\beta\right\rfloor }{\left\lceil \frac{g^{\beta}}{\ln g} \right\rceil}  
			\\&\leq \left({e\ln g} \right)^{\frac{g^\beta}{\ln g}} \cdot (2\pi \frac{g^\beta}{\ln g})^{-1/2} \cdot \exp\left(- \frac{g^\beta}{2\ln^2 g}(1 + o(1))\right)
			\\&= O\left(\exp\left({2\frac{g^\beta}{\ln g}\ln\ln g}\right)\right),
		\end{align*}
		and so \begin{align*}\left(\star\right) &=  O\left(\exp\left({2\frac{g^\beta}{\ln g}\ln\ln g}\right)\right) \exp\left(-\frac{2\alpha^3\zeta_g^2g^\beta}{e^{C\Delta_+}(1-2\beta) \ln g} + o\left(\frac{\alpha^3\zeta_g^2g^\beta}{e^{C\Delta_+}\ln g}\right)\right)\\ &= 
			O\left(\exp\left(\frac{2g^\beta}{\ln g}\left(\ln\ln g-\frac{\alpha^3\zeta_g^2}{e^{C\Delta_+}(1-2\beta)} + o\left(\frac{\alpha^3\zeta_g^2}{e^{C\Delta_+}}\right)\right)\right)\right)\\& = o\left (g^{-3}\right)
		\end{align*}
		because \textbf{$\alpha^3\zeta_g^2 \gg e^{C\Delta_+} \ln \ln g$}.
	\end{proof}
	
	\subsubsection{Tools for the lower bound}
	The aim of this section is to prove Lemma \ref{lem:detmin}, which estimates from below the distances from $\rho$ to the $g^{\frac{1}{2}-\frac{1}{\ln^{3/4}g}}$ first vertices discovered during the exploration around $\rho$. As for Lemma \ref{lem:detmaj}, we start with a coupling lemma.
	
	\begin{lemma}\label{lem:couplage2}
		Let $(G, w)$ be the surface obtained after finding $g^{\frac{1}{2}-\frac{1}{\ln^{3/4}g}}$ vertices in the exploration. Suppose that the exploration around $\rho$ does not stop prematurely, and that fewer than $b$ bad steps are made before we entirely discover $G$.
		Then, there is an injection $\iota\colon \left(G, \rho\right) \hookrightarrow \left(\widehat S, \rho\right)$, where $\widehat S$ is a copy of $\widehat{\mathfrak S}_\nu$, such that \[\forall x \in G, \quad d_+^{G}(x) \geq d^{S}(\rho, \iota(x))- 5\Delta_+ b.\] 
	\end{lemma}
	
\begin{figure}[ht]
		\centering
		\subfloat[][When $x_0 \notin P'$, we concatenate the blue path with parts of the red one. The boundary corresponding to the bad step is in bold, green.]
		{			\centering
			\includegraphics[width=0.35\textwidth]{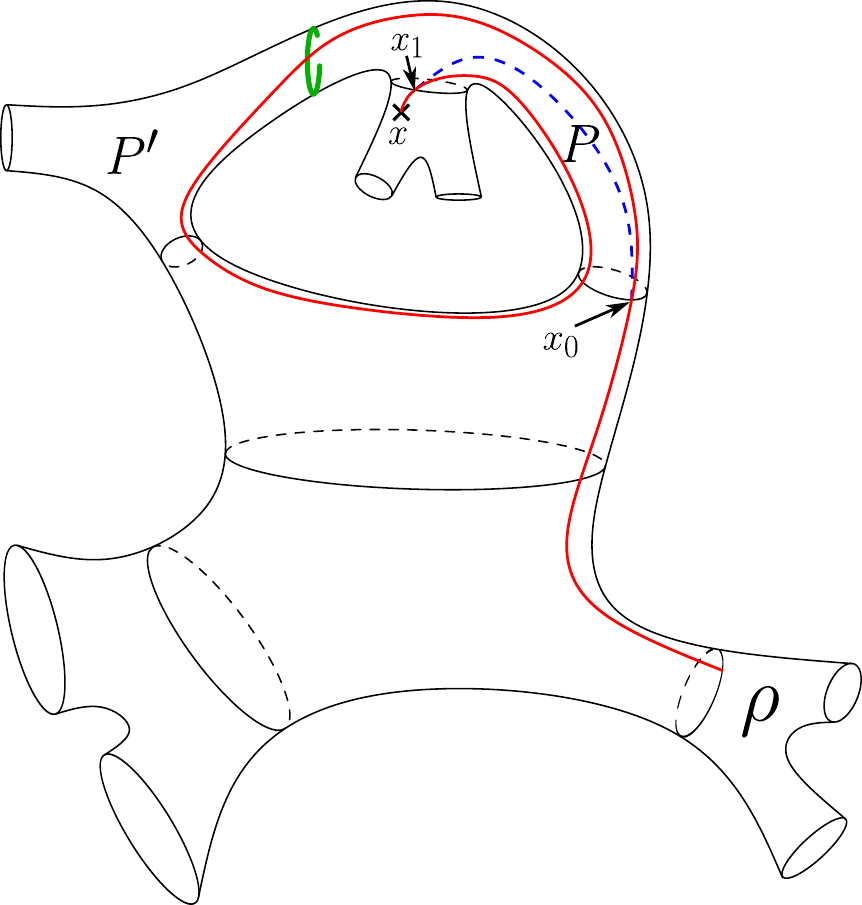}
			\label{fig:path2}
		}
		\hfill
		\subfloat[width = 0.3\textwidth][When $x_0 \in P'$, we replace the red path with the blue one until we leave $P$. The boundary corresponding to the bad step is in bold, green.]
		{
			\centering
			\includegraphics[width=0.35\textwidth]{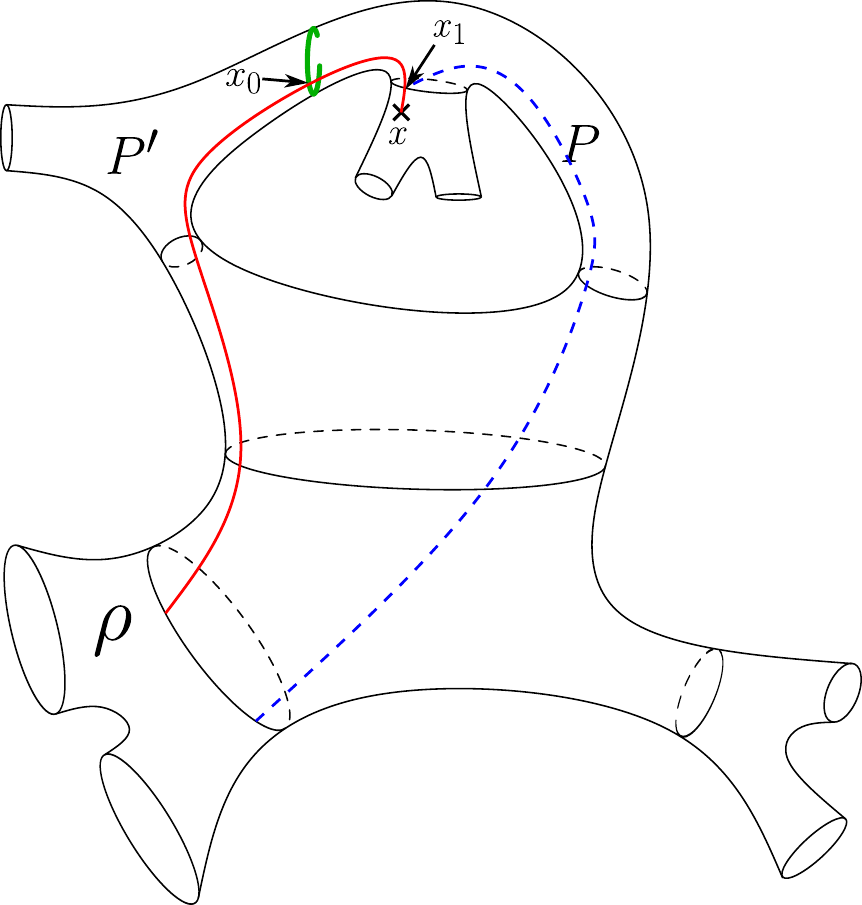}
			\label{fig:path1}
		}
		\caption{The two ways of modifying $\gamma$ in the proof of Lemma \ref{lem:couplage2}.}
		\label{fig:path1and2}
	\end{figure}

	\begin{proof}
		$\widehat S$ is built by "forgetting" the bad steps in the exploration. Start from a copy of $\rho$. Suppose that when doing the exploration of $G$, we want to pair a half-edge of $P$ to another pair of pants $P'$. If it is not a bad step, then glue a copy of $P'$ to $\iota(P)$  in $\widehat S$ (using the same length and twist). On the contrary, if it corresponds to a bad step, \textit{i.e.} $P'$ was already discovered previously in the exploration, then instead of gluing the corresponding two pairs of pants together in $\widehat S$ as we did in $G$, we glue to each of them an independent copy of $\mathfrak S_\nu$ and redraw an independent weight for both half-edges (\textit{i.e.} for both boundary components, which are not glued together). The image $\iota(x)$ of $x \in P\cap P'$ is chosen to be the corresponding point in $\iota(P)$. We complete this subsurface to obtain a copy of $\widehat{\mathfrak S}_\nu$.
		
		Take $x \in G$, and suppose that $d_+^{G}(x) = R$. Let us now consider a path $\gamma$ from $\rho$ to $x \in G$ whose maximum possible length is $R$. Our aim is to build a path from $\iota(\rho)$ to $\iota(x)$ of length less than $R+5\Delta_+ b$.
		
		First, we build another path $\psi$, still from $\rho$ to $x$, of length at most $R + 3\Delta_+b$, which does not use bad steps (it never passes through a gluing due to a bad step).
		More precisely, we do it by induction, proving that if $\gamma$ uses only the $k$ first bad steps and is of length at most $R$, then we can find another path $\widetilde{\psi}$ using only the $k-1$ first bad steps of length at most $R+3\Delta_+$. Then we obtain $\psi$ recursively, because at most $b$ bad steps are made.
		
		Suppose that the $k$-th bad step glues together $P$ and $P'$. Let $x_1$ the last point of $\gamma$ in $\left(P \cup P'\right)$, $x_1 = \gamma(t_1)$. Without loss of generality, we suppose $x_1 \in P$. Let $x_0$ be the first point of $\gamma$ in $P$, $x_0 = \gamma(t_0)$. \begin{itemize}
			\item if $x_0 \notin P'$, we take a path $\widetilde{\gamma}$ inside $P$ from $x_0$ to $x_1$ of length less than $\Delta_+$. Concatenate $\gamma\vert_{[0, t_0]}$, $\widetilde{\gamma}$ and $\gamma\vert_{[t_1, 1]}$ to obtain $\widetilde{\psi}$ (see Figure \ref{fig:path2}).
			\item if $x_0 \in P'$, let us write $r$ for the maximum possible length of $\gamma\vert_{[0, t_0]}$. By definition of a bad step, there exists another boundary of $P$ which has been paired before $P \cap P'$. According to the second part of Lemma \ref{lem:breadth}, it implies that there is a point $x_0'$ on this boundary at distance at most $R+2\Delta_+$ from $\rho$. Concatenate any path from $\rho$ to $x_0'$, a path from $x_0'$ to $x_1$, and $\gamma\vert_{[t_1, 1]}$ to obtain $\widetilde{\psi}$ (see Figure \ref{fig:path1}). 
		\end{itemize}
		
		The projection $\iota(\psi)$ of $\psi$ in $\widehat S$ is a continuous path from $\rho$ to $\iota(x)$ in $\widehat S$. The length of the part of $\iota(\psi)$ in pairs of pants without bad steps is smaller than the corresponding parts of $\psi$, so is at most $R+3\Delta_+b$. We modify $\iota(\psi)$ a last time to ensures for any pair of pants $\iota(P)$ affected by a bad step, the length of $\iota(P) \cap \iota(\psi)$ is at most $\Delta_+$ -- if it not the case, create a shortcut between the first and the last point of $\iota(P) \cap \iota(\psi)$ As a consequence, because at most $b$ bad steps are made and at most two pairs of pants are affected by each one of them, this last path is of length at most $R + 5\Delta_+b$, as requested.
	\end{proof}
	
	\begin{proof}[Proof of Lemma \ref{lem:detmin}]
		We use the coupling Lemma \ref{lem:couplage2} to pair the subsurface $G$ obtained at the end of the partial exploration with a copy $\widehat S$ of $\widehat{\mathfrak S}_\nu$.
		
		Let $R$ be the maximum possible maximum distance from $\rho$ reached in $G$, \textit{i.e.} \[R = \sup_{x \in G} d_+^{G}(x).\] 
		As $g^{\frac{1}{2}-\frac{1}{\ln^{3/4}g}}$ vertices are at distance less than $R$ from $\rho$ in $S_1$, and because less than $\ln^{3/4} g$ bad steps are made, we should have that at least $g^{\frac{1}{2}-\frac{1}{\ln^{3/4}g}}$ pairs of pants are at distance less than $R+5\ln^{3/4}g\Delta_+$ from $\rho$ in $\widehat S$, \textit{i.e.} $\widehat{N}_{R +3\Delta_+\ln^{3/4} g} \geq g^{\frac{1}{2}-\frac{1}{\ln^{3/4}g}}$ in $\widehat S$.
		
		Take $\widetilde{R} = \frac{1}{2\alpha}\ln g  - \zeta_g$, and $\varepsilon = \left(\frac{1}{2\widetilde{R}}-\frac{1}{\widetilde{R}\ln^{3/4} g}\right)\left(\ln g - \alpha\right) = \frac{2\alpha^{2} \zeta_g}{\ln g} + o\left(\frac{\alpha^{2} \zeta_g}{\ln g}\right),$ so that $0<\varepsilon<1$, $\varepsilon \limitginf{}{} 0$ because \textbf{$\ln^{1/4} g \ll \alpha\zeta_g\ll \ln g$}.
		As in Lemma \ref{lem:detmaj}, since $\zeta_g \gg \frac{e^{C\Delta_+}}{\alpha^{2}}\ln^{3/4}g$, Theorem \ref{th:concentration} and Corollary \ref{coro:surfacecomplete} apply and so	
		\begin{align*}\p\left (\widehat{N}_{\widetilde{R}} \geq g^{\frac{1}{2}-\frac{1}{\ln^{3/4}g}} \right)
			&\leq 2\exp\left(-\frac{\alpha^3\zeta_g^2}{e^{C\Delta_+}\ln g} + o\left(\frac{\delta_-\alpha^3\zeta_g^2}{ \Delta_+\ln g}\right)\right) = o\left(1\right),\end{align*}
		because $\alpha^3\zeta_g^2\gg e^{C\Delta_+}\ln g$.
		So $R \geq \frac{1}{2\alpha}\ln g  - \zeta_g$ with probability $1-o(1)$. According to Lemma \ref{lem:breadth}, any unpaired half-edge is at least $R-2\Delta_+$ away from $\rho$ in $\mathfrak S_{\nu, g}$. We get that with probability $1-o\left(1\right)$, any unpaired boundary component after the exploration of $g^{\frac{1}{2}-\frac{1}{\ln^{3/4}g}}$ vertices is at least at distance \[\frac{1}{2\alpha}\ln g - \zeta_g - \Delta_+\left(2+5\ln^{3/4} g\right) \geq \frac{1}{2\alpha}\ln g - 2\zeta_g,\] from $\rho$.
	\end{proof}
	
	\subsection{Proof of the main theorem}\label{Prooffin}
	
	The aim of this section is to give the proof of Theorem \ref{th:Simp}. The main idea, already used in \cite{Bollobas} and \cite{BCP_mindiam}, is to observed that by a birthday paradox kind of argument, if we do two exploration around two different pairs of pants P and $P'$, we should wait until around $\sqrt{g}$ vertices have been found in which to see them merge. Lemmas \ref{lem:detmaj} and \ref{lem:detmin} are then used to estimate the radius of the neighbourhood discovered during these explorations.
	
	\begin{proof}[Proof of Theorem \ref{th:Simp}]
		The proof is divided into two parts, in which we separately bound from above and from below the diameter of $\mathfrak S_{\nu, g}$. Each part is based on one of the Lemmas \ref{lem:detmaj} or \ref{lem:detmin} proved in Section \ref{Algo}.
		\begin{itemize}
			\item First, we show the upper bound \[\diam \mathfrak S_{\nu, g} \leq \frac{1}{\alpha} \ln g + 2\zeta_g.\]
			
			Take two pairs of pants $P, {P}'$. We proceed with the exploration around $P$ as described in Definition \ref{defi:algo}, stopped when $\sqrt{g}\ln g$ vertices are discovered. Let us suppose that the surface is not disconnected. Then, with probability $1-o\left(g^{-2}\right)$, at least $\sqrt{g}\ln g - 3 \ln^3g > \frac{1}{2}\sqrt{g}\ln g$ half-edges are discovered but not paired yet. Indeed, at each step, the number of half-edges discovered but not paired yet increases by 1 (gains 2, loses 1) unless it was a bad step (gains 0, loses 2), and by Lemma \ref{lem:borneserreurs}, at most $\ln^3 g$ bad steps are made.
			Let us now proceed with a similar exploration around $P'$. At each step, the probability that it does not merge with the exploration around $P$ is bounded from above by \[1-\frac{\sqrt{g}\ln g}{2g}.\]
			As a consequence, the probability that we see $\sqrt{g}\ln g$ vertices around $P'$ but it does not merges with the exploration around $P$ is bounded from above by \[\left(1-\frac{\ln g}{2\sqrt{g}}\right)^{\sqrt{g}\ln g} \leq \exp\left(-\frac{\ln^2 g}{2}\right) = o\left(g^{-2}\right).\]
			
			Therefore, with probability $1-o\left(g^{-2}\right)$, either the surface is disconnected, or the two explorations merge. Let us now use together Lemmas \ref{lem:borneserreurs} and \ref{lem:detmaj}. They implies that with probability $1-o(g^2)$, under the event that the surface is connected, any of the half-edges we paired during the two explorations is at distance less than $\frac{1}{2\alpha} \ln g + \zeta_g$ from the starting point $P$ or $P'$ of the associated exploration. In particular, when the surface is connected, with probability $1-o\left(g^{-2}\right)$, \[d_h(P, \widetilde P) \leq \frac{1}{\alpha} \ln g + 2\zeta_g.\] By a union bound on the pairs of pants $(P, P')$, with probability $1-o(1)$, either the surface is connected, or its diameter is bounded from above by \[\frac{1}{\alpha} \ln g + 2\zeta_g.\] Because the graph is connected with high probability (see, for example, \cite[Corollary 1.1]{WORMALD1981156}), we get the upper bound.
			\item Next, we show the lower bound \[\diam \mathfrak S_{\nu, g} \geq \frac{1}{\alpha} \ln g - 2\zeta_g.\]
			As for the upper bound, consider again the exploration around $P$, but this time, we stop it as soon as $g^{\frac{1}{2}-\frac{1}{\ln^{3/4} g}}$ vertices are found. It take less than $3g^{\frac{1}{2}-\frac{1}{\ln^{3/4} g}}$ steps, and so after that, at most $3g^{\frac{1}{2}-\frac{1}{\ln^{3/4} g}}$ half-edges are discovered but not paired yet. Proceed similarly with some $P'$ which is not in the explored part. At each step, the probability that the current half-edge is paired with one of the first exploration is bounded from above by \[ C\frac{g^{1/2-\frac{1}{\ln^{3/4} g}}}{g}\,\] for some $C>0$ independent of $g$, as in the proof of Lemma \ref{lem:borneserreurs}. Consequently, the probability that the two explorations merge is bounded from above by \[ \left(C\frac{g^{1/2-\frac{1}{\ln^{3/4} g}}}{g} g^{1/2-\frac{1}{\ln^{3/4} g}}\right) = Cg^{-\frac{2}{\ln^{3/4} g}} = o(1).\]
			
			According to Lemma \ref{lem:detmin}, with probability $1-o(1)$, these two balls contain all the points at distance $\frac{1}{2\alpha}\ln g - \zeta_g$ of $P$ (resp. $P'$). In particular, with probability $1-o(1)$, the distance between $P$ and $P'$ is greater than $\frac{1}{\alpha}\ln g - 2\zeta_g.$
		\end{itemize}
		As a consequence, with high probability, \[\frac{1}{\alpha}\ln g - 2\zeta_g \leq \diam \mathfrak S_{\nu, g} \leq \frac{1}{\alpha}\ln g + 2\zeta_g,\] hence the result.
	\end{proof}

	\bibliographystyle{alpha}
	\bibliography{ref.bib}
	%\printbibliography
	
	\newpage
	\appendix
	
	\section{Reminders of probability}
	\label{Probrem}
	
	In this section, we recall the two concentration inequalities used in the paper.
	\begin{proposition}[{McDiarmid's inequality}]\label{prop:McD}
		Let $X_1, \ldots, X_n$ be independent random variables, with $X_i$ taking its values in $\mathcal X_i$. We fix $f \colon \mathcal X_1 \times \ldots \times \mathcal X_n \to \R$, such that \[\forall i, \exists c_i \in \R,\ \forall x_1, \ldots, x_n, \forall x_i',\quad \left\lvert f(x_1, \ldots x_i, \ldots, x_n) - f(x_1, \ldots, x_i', \ldots x_n)\right\rvert \leq c_i.\]
		Then \[\p\left(\left\lvert f(X_1, \ldots, X_n) - \E f(X_1, \ldots, X_n)\right\rvert \leq \varepsilon\right) \leq 2\exp \left(-\frac{2\varepsilon^2}{\sum_{i = 1}^{n} c_i^2} \right).\]
	\end{proposition}
	
	\begin{proposition}[{Paley-Zigmund's inequality}]\label{prop:PZ}
		Let $Z$ be a non-negative random variable. Then for all $0 \leq \theta \leq 1$, \[\p \left(Z \geq \theta\ \E Z \right) \geq (1-\theta)^2\, \frac{\E\left[Z\right]^2}{\E Z^2}.\]
	\end{proposition}

	\section{Some facts on hyperbolic surfaces}
	\label{hypreminder}
	
	In this appendix, we collect some "basic" facts about hyperbolic surfaces. The first subsection is devoted to the collar lemma and how it is used to construct surfaces of large diameter. The second deals with the proof of Proposition \ref{prop:Ddelta}, about the geometry of pairs of pants, and more especially about estimates on the diameter of pairs of pants and on the minimum distance between two boundary components.
	
	\subsection{Collar lemma}
	\label{collar}
	\begin{lemma}[Collar Lemma -- \cite{Buser}, Theorem 4.1.1]\label{lem:collar}
		Let $S$ be a closed hyperbolic surface, and $\gamma$ a simple closed geodesic on $S$. Denote $\ell(\gamma)$ the length of $\gamma$, and $w(\gamma) = \arcsinh\left(\frac{1}{\sinh\left(\frac{1}{2}\ell(\gamma)\right)}\right)$. Then, the set \[\left\lbrace p\in S \middle \vert d(p, \gamma) < w(\gamma) \right\rbrace\] is isometric to the cylinder $[-w(\gamma), w(\gamma)] \times \mathbb S_1$ with the Riemannian metric \[\mathrm{d}s^2 = \mathrm{d}\rho^2 + \ell^2(\gamma)\cosh^2(\rho)\, \mathrm{d}dt^2.\]
	\end{lemma}
	
\begin{figure}[ht]
		\centering
		\includegraphics[width=0.8\textwidth]{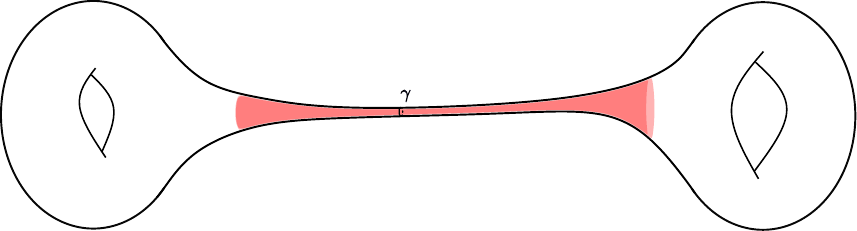}
		\caption{Coloured in red, the cylinder associated to the curve $\gamma$ in Lemma \ref{lem:collar}. Because $\gamma$ is "short", the width of the cylinder is large, and so is the diameter of the whole surface. In general, the same works for surfaces of any genus as soon as $\gamma$ disconnects the surface.}
		\label{fig:largediameter}
	\end{figure}

	To build surfaces with a large diameter, we can glue pairs of pants with very small boundary components. Indeed, this lemma tells us that when $\ell(\gamma) \to 0$, on the contrary, the width of the associated cylinder $w(\gamma) \to \infty$. For example, a surface of genus 2 with a large diameter could look like Figure \ref{fig:largediameter}.
	
	\subsection{Pair of pants}
	\label{pantalon}
	As we see in Sections \ref{TLS} and \ref{Surface}, the geometry of pairs of pants is the base of our proof. We recall that a pair of pants can be built by gluing together two identical right-angled hyperbolic hexagons, and so its geometry is well known. For general facts on pairs of pants and right-angled hexagons, we refer to \cite[Chapters 2-3]{Buser}. In this appendix, we only prove the estimates for the diameter $\Delta_P$ and the minimum distance $\delta_P$ between two boundary components of a pair of pants $P$ given in Proposition \ref{prop:Ddelta}. Note that when the boundary length is fixed to a value $l>0$, we use the notation $\delta_l = \delta_+ = \delta_-$.
	
	\begin{proof}[Proof of Proposition \ref{prop:Ddelta}]
		The expressions of $\delta_-, \delta_+$ come from the fact that for a pair of pants with boundary lengths $2l_1, 2l_2, 2l_3$, the distance between the boundary components $1$ and $2$ is given by \[\cosh^{-1}\left(\frac{\cosh l_3 + \cosh l_1 \cosh l_2}{\sinh l_1 \sinh l_2} \right).\] We minimize (respectively maximize) this quantity to obtain $\delta_-$ (respectively $\delta_+$). From this, we have that \[\cosh(\delta_-) \geq \frac{1+2e^{l^- - 2l^+}}{1 - 2e^{-2l^+}+e^{-4l^+}}\geq 1+2e^{l^- - 2l^+} \in [1, 3].\] Using concavity of $\cosh^{-1}$ between $1$ and $3$, we get that \[\delta_- \geq  e^{l^- - 2l^+}\cosh(3) \geq e^{- 2l^+}.\]
		On the other hand, because \[\forall x \in \Rp, \quad \cosh^{-1}(x) \leq \ln 2x, \quad \cosh x \leq e^x\quad \text{and} \quad \sinh^2x \geq x^2,\] we obtain immediately that \[\delta_+ \leq \ln\frac{2\left (e^{l^+} + e^{2l^-}\right )}{\left(l^-\right)^2} \leq \ln{\frac{4e^{2l^+}}{\left(l^-\right)^2}}.\]
		The lower bound on $\Delta_+$ is immediate. For the upper bound, remark that the diameter of a pair of pants is bounded by twice the diameter of one of the two right angled hyperbolic hexagons used to build it, which is itself bounded by $2l^+ +2\delta_+\leq 4\max(l_+, \delta_+)$.
		The universal lower bound on $\Delta_+$ comes from the fact that it is at least equal to the diameter of a pair of pants all of whose boundary components are of length $l$ for any $l \in [2l^-, 2l^+]$. In such a case, we recall that \[\delta_l = \cosh^{-1}\left(\frac{\cosh l + \cosh^2 l}{\sinh^2 l} \right).\]
		Comparing $\frac{\cosh\delta_l}{\cosh{l}}$ to $1$, we get that $\delta_l = l$ for $l = \ln(2+\sqrt{3})$. Because $\delta_l$ is a decreasing function of $l$, this also gives a universal lower bound for $\max(l, \delta_l)$ and so for $\Delta_+$.
		Finally, \begin{align*}\ln\left(\frac{\Delta_+}{\delta_-} +1\right) &\leq \ln\left(\frac{1}{\delta_-}\right) + \ln\left({\Delta_+} +{\delta_-}\right) \\&\leq 2l^+ + \ln\left(2{\Delta_+}\right)\\ &\leq 4\Delta_+, 
		\end{align*} which concludes.
	\end{proof}
	
	\section{Geometric filtration and spatial Markov property}
	\label{geometricmarkov}
	In this appendix, we prove a geometric Markov property for our model, using the notions introduced in Section \ref{TLS}, especially in Definition \ref{defi:sttr} (stopping trees and associated $\sigma$-algebra). We recall that a complete subtree of $(\mathfrak B, \rho)$ is a subtree (rooted at $\rho$) such that each node is either a leaf or of degree 2. For such a finite complete subtree $T$, we recall the definition of $^{\mathtt C}{}{T}$ and $\mathfrak B \setminus ^{\mathtt C}{}{T}$ given in Definition \ref{defi:childtree}. Then, $\mathfrak B \setminus ^{\mathtt C}{}{T}$ is a forest of copies of $\mathfrak B$ rooted at the grandchildren of the leaves of $T$. The twists and length determining $T$ (\textit{i.e.} in $\mathcal F_T$), and the ones in each of these copies are independent by construction -- this is exactly why we want to consider $^{\mathtt C}{}{T}$ instead of $T$.
	
	From this setting we can directly obtain a weak Markov property, and deduce a strong one.
	
	\begin{theorem}[{Weak Markov property}]\label{th:Markovfaible}
		Let $T \preceq \mathfrak B$ be a (deterministic) finite complete subtree, with $L = L(T)$ leaves. Then $\mathfrak B \setminus ^{\mathtt C}{}{T}$ can be written as $4L$ iid copies of $\mathfrak B$, denoted $\mathfrak B_1^{T}, \ldots \mathfrak B_{4L}^{T}$. In addition, for any family $(f_i)_{i \in \N}$ of positive functions, 
		\[\E\left[\bigop{\prod}{i=1}{4L} f_i\left(\tau^{\mathfrak B_i^{T}}\right)\ \middle|\ \mathcal F_T\right] = \bigop{\prod}{i=1}{4L} \E f_i\left(\tau^{\mathfrak B}\right).\]
	\end{theorem}

	\begin{proof}
		It is just a reformulation of the fact that the trees are all identically distributed. The completeness of the tree guarantees that all leaves, and \textit{a fortiori} their grandchildren, form an anti-chain, so that the associated subtrees are independent from one another. 
	\end{proof}

	\begin{theorem}[{Strong Markov property}]\label{th:Markovfort}
		Let $T$ be a complete finite stopping tree. Then the number $L(T)$ of leaves of $T$ is $\mathcal F_T$-measurable, and  $\mathfrak B \setminus ^{\mathtt C}{}{T}$ can be written as a disjoint union of $4L$ independent copies of $\mathfrak B$, denoted $\mathfrak B_1^{T}, \ldots \mathfrak B_{4L}^{T}$.
		
		In addition, for any family $(f_i)_{i \in \N}$ of positive functions, 
		\[\E\left[\bigop{\prod}{i=1}{4L(T)} f_i\left(\tau^{\mathfrak B_i^{T}}\right)\ \middle|\ \mathcal F_T\right] = \bigop{\prod}{i=1}{4L(T)} Ef_i\left(\tau^{\mathfrak B}\right).\]
	\end{theorem}
	
	\begin{proof}
		For any deterministic $(T_0, \rho) \preceq (\mathcal B, \rho)$, the event $\lbrace T = T_0\rbrace$ is $\mathcal F_T$-measurable. Consequently, so is $L(T) = \sum_{T_0} L(T_0)\mathds{1}_{T=T_0}$.
		
		In this way, the random variable $\bigop{\prod}{i=1}{4L(T)} \E f_i\left(\tau^{\mathfrak B_i^T}\right)$ is $\mathcal F_T$-measurable. It suffices to check the characteristic property of conditional expectation to conclude. Take $A \in \mathcal F_T$.
		
		\begin{align*}
			\E\left[\mathds{1}_A \bigop{\prod}{i=1}{4L} f_i\left(\tau^{\mathfrak B_i^{T}}\right)\right] &= \bigop{\sum}{{\#T_0 < \infty}}{} \E\left[\mathds{1}_{A \cap \{T = T_0\}} \bigop{\prod}{i=1}{4L(T_0)} f_i\left(\tau^{\mathfrak B_i^{T_0}}\right)\right]\\ &= \lim_{n\to\infty} \bigop{\sum}{\#T_0 \leq n}{} \E\left[\mathds{1}_{A \cap \{T = T_0\}} \bigop{\prod}{i=1}{4L(T_0)} f_i\left(\tau^{\mathfrak B_i^{T_0}}\right)\right]\\&= \lim_{n\to\infty} \bigop{\sum}{\#T_0 \leq n}{} \E\left[\mathds{1}_{A \cap \{T = T_0\}}\ \E\left[\bigop{\prod}{i=1}{4L(T_0)} f_i\left(\tau^{\mathfrak B_i^{T_0}}\right)\ \middle|\ \mathcal F_{T_0}\right]\right],\end{align*}
		using that $\mathds{1}_{A \cap {T = T_0}}$ is $\mathcal F_{T_0}$-measurable, because $T$ is a stopping-tree. Thanks to the weak Markov Property \ref{th:Markovfaible}, we get \begin{align*}
			\E\left[\mathds{1}_A \bigop{\prod}{i=1}{4L} f_i\left(\tau^{\mathfrak B_i^{T}}\right)\right] &= \lim_{n\to\infty} \bigop{\sum}{\#T_0 \leq n}{} \E\mathds{1}_{A \cap \{T = T_0\}} \bigop{\prod}{i=1}{4L(T_0)} \E f_i\left(\tau^{\mathfrak B}\right)\\ &= \E\left[1_A \bigop{\prod}{i=1}{4L(T)} \E f_i\left(\tau^{\mathfrak B}\right)\right].
		\end{align*}
	\end{proof}
	
	To finish this part, we will show that the trees $T_R$, $T'_R$ defined after Definition \ref{defi:sttr} are stopping trees.

	\begin{proposition}\label{prop:stoptime}
		$T_R$, $T'_R$ are stopping trees.
	\end{proposition}
	
	\begin{proof}
		Let us fix some deterministic subtree $T$ of $\mathfrak B$. Then $T_R$ is not contained in $T$ means that $U_R$ is not contained in $T$, and so that there exists a leaf $P$ of $T$ such that \[d_h(\partial \rho ^-, \partial P^+) < R.\] The last event is $\mathcal F_{T}$-measurable, and so is $\mathds{1}_{T_R \preceq T}$. We do a similar reasoning for $T'_R$, using the fact that $U'_R$ is entirely determined by $\mathcal S_R$, which can exactly be described in terms of the metric on $\mathcal B_R \preceq T'_R$.
	\end{proof}
	
\end{document}